\documentclass[twoside,11pt]{article}

\usepackage[abbrvbib]{jmlr2e}

\usepackage{jmlr2e}

\usepackage{amsmath}

\usepackage{tikz}
\usetikzlibrary {shapes}
\usetikzlibrary {arrows}
\usetikzlibrary {positioning}
\usetikzlibrary {calc}
\usetikzlibrary{fit}					
\usetikzlibrary{backgrounds}

\usepackage{color}

\newcommand{\pa}{\mathrm{pa}}

\newcommand{\an}{\mathrm{an}}

\newcommand{\pst}{\mathrm{pst}}
\newcommand{\mb}{\mathrm{mb}}
\newcommand{\dis}{\mathrm{dis}}

\newcommand{\dse}{\,\mbox{$\perp$}\,}

\newcommand{\cip}{\mbox{\,$\perp\!\!\!\perp$\,}}
\newcommand{\sk}{\mathrm{sk}}

\newcommand{\cd}{\,|\,}

\newcommand{\ci}{\mbox{\protect $\: \perp \hspace{-2.3ex}
\perp$ }}

\newcommand{\notdse}{\nolinebreak{\not\hspace{-1.5mm}\dse}}

\newcommand{\notci}{\nolinebreak{\not\hspace{-1.5mm}\ci}}

\newcommand{\n}[0]{\hspace*{.35em}}

\newcommand{\nn}[0]{\hspace*{.7em}}

\newcommand{\node}{\mbox {\LARGE
{$\mbox{$\circ$}$}}}

\newcommand{\ful}{\mbox{$\, \frac{ \nn \nn \;}{ \nn \nn
}$}}

\newcommand{\fla}{\mbox{$\hspace{.05em} \prec
\!\!\!\!\!\frac{\nn \nn}{\nn}$}}

\newcommand{\fra}{\mbox{$\hspace{.05em} \frac{\nn
\nn}{\nn
}\!\!\!\!\! \succ \! \hspace{.25ex}$}}

\newcommand{\arc}{\mbox{$\hspace{.06em} \prec
\!\!\!\!\!\frac{\nn \nn}{\nn}
\!\!\!\!\!
\succ\! \hspace{.25ex}$}}

\newtheorem{eexample}[theorem]{Example}

\usepackage{enumerate}

\usepackage{lastpage}
\jmlrheading{23}{2022}{1-34}{4/21; Revised 9/21}{2/22}{21-0425}{Kayvan Sadeghi and Terry Soo}

\ShortHeadings{Conditions for Causal Structure Learning}{Sadeghi and Soo}
\firstpageno{1}

\begin{document}

\title{Conditions and Assumptions for Constraint-based Causal Structure Learning}

\author{\name Kayvan Sadeghi \email k.sadeghi@ucl.ac.uk \\
       \addr Department of Statistical Science\\
       University College London\\
       London, United Kingdom
       \AND
       \name Terry Soo \email math@terrysoo.com \\
       \addr Department of Statistical Science\\
       University College London\\
       London, United Kingdom}

\editor{Peter Spirtes}

\maketitle
\begin{abstract}
 We  formalize constraint-based structure learning of the ``true" causal graph from observed data when unobserved variables are also existent. We provide conditions for a ``natural" family of constraint-based structure-learning algorithms that output graphs that are Markov equivalent to the causal graph. Under the faithfulness assumption, this natural family contains all exact structure-learning algorithms.  We also provide a set of assumptions, under which any natural structure-learning algorithm outputs Markov equivalent graphs to the causal graph.
 These assumptions can be thought of as a relaxation of faithfulness, and most of them can be directly tested from (the underlying distribution) of the data, particularly when one focuses on structural causal models. We specialize the definitions and results for structural causal models.
\end{abstract}


\section{Introduction}
Inferring causal relationships has always been one of the main objectives in many fields of study, ranging from philosophy to  science \citep{MR867618}. In the past decades, \emph{causal models} \citep{pea09} have been used to infer causal relationships from observed data. Extensive research has been conducted on defining, interpreting, and applying causal models, and these models have become
mainstream in statistics and computer science \citep{pea88,spioo,daw00,hit18}. Today, a very popular method for inferring causal relationships is based on the use of statistical models over graphs with nodes that are random variables representing the quantities of interest. Edges indicate probabilistic  dependence among these variables conditional on some other variables in the graph, which with some additional assumptions can be interpreted as causal relationships. These models are called \emph{graphical (Markov) models} \citep{lau96}.

\emph{Structural causal models} (SCMs), also known as \emph{structural equation models},  are one of the most used causal models \citep{spioo,pea88}. SCMs have been generally defined for directed acyclic graphs (DAGs) (i.e., a Bayesian network), where an arrow implies a ``direct"  cause. However, for the purpose of causal inference, it may not be possible to determine or measure all  common causes.
For this reason, SCMs associated to DAG models with latent variables (represented by ancestral graphs \citep{ric02} without undirected edges) have been defined \citep{spioo,zha08b}.

A main goal of causal inference, based on graphical models, is causal discovery; in this framework we are concerned with finding graphs, using only observed data,  that represent the causal structure, which is represented by an unknown ``true" causal graph.    Algorithms for constructing or selecting such a graph are known as \emph{structural learning} algorithms.  There has been a surge in defining different structural learning algorithms; they are generally categorized as constraint-based \citep{spioo}, score-based \citep{chi03} --which work on the factorization of density--, or hybrid \citep{tsa06} --which is a combination of these.

We are concerned with  \emph{constraint-based} structure-learning algorithms for the general case of ancestral graphs without undirected edges. These algorithms use conditional independence statements, which are
theoretically
testable from observed data to find graphs that, at best, are Markov equivalent to the true causal graph.   A main sufficient condition for the algorithms to be \emph{correct} or  \emph{work}, that is, to provide an output that is Markov equivalent to the causal graph, is faithfulness  \citep{zha08,woo98,ste06,col14}; this sufficient condition relates the  independencies induced by the distribution (of the SCM) and the separations induced by the causal graph.   However,  faithfulness to the unknown causal graph, or similar  assumptions \citep{zha03}, are generally not testable. In addition, the (sufficient) assumption of faithfulness is in some circumstances strong, and can be relaxed;  see Section \ref{sec:alg}.

In this paper, we formalize the problem of constraint-based structure learning from observed data.   Most constraint-based structure-learning algorithms contain two tasks which may be performed in parallel or separately:  finding the undirected edges of the graph, the skeleton;   and directing the edges of the skeleton. We show that under the assumption of faithfulness and using a formalization of faithfulness,  provided in \citet{sad17}, if the output of such algorithms works, then pair of  conditions known as \emph{ordered stability} must be satisfied with respect to the output;   see Section \ref{sec:ouds} and  Remark \ref{must}.   Ordered stability is a condition that informs us how the edges of the graph should be directed and we propose the class of \emph{natural} structure-learning algorithms, which is one that should respect stability when directing the edges; see  Definition \ref{algorithm}. We stress that in this paper, we do not provide specific algorithms that could be implemented, but concentrate on providing conditions that structure-learning algorithms should satisfy.  More importantly, we provide clear assumptions, many of which
are directly testable from the distribution, under which a natural structure-learning algorithm generates the correct output.

We provide the theory for the general class of models under the assumption that the distribution is Markovian to the true causal graph; see Theorem \ref{coro:main}. The assumptions are the following: the distribution satisfies ordered stability and the converse of the pairwise Markov property with respect to the true causal graph, and also the distribution satisfies what we call path-stability.

We will make precise our assumptions later, and we note that they can be thought of as a relaxation of faithfulness, and in particular, in the case of DAGs our assumptions are logically implied by faithfulness; see Corollary \ref{logically-implied}.  We also specialize the definitions and results for SCMs; see Theorem \ref{main-scm-cor}, and, in particular, provide testable sufficient conditions for the converse of the pairwise Markov property to be satisfied.
We hope that the framework provided here will be a useful theoretical tool in defining appropriate SCMs and  concrete algorithms with strong theoretical backing.

%

The structure of the paper is as follows: In the next section, we provide known definitions and results needed for the theory presented in this paper. In Section \ref{sec:gen-set}, we formalize the settings of the problem. Section \ref{sec:faith} deals with the theory under the faithfulness assumption. In Section \ref{sec:alg}, we introduce the structure learning theory without the faithfulness assumption. Section \ref{sec:st-l-SCM} specializes on causal structure learning for SCMs. We end the paper with a summary and discussion in Section \ref{sec:summary}.
\section{Preliminaries}
In this section, we provide the basic definitions and concepts as well as the known results needed for the paper.
\subsection{Graph theoretical definitions and concepts}
We usually refer to a graph as an ordered pair $G=(V,E)$, where $V$ is the \emph{node} set and $E$ is the \emph{edge} set. When nodes $i$ and $j$ are the endpoints of an edge, we call them \emph{adjacent}, and write $i\sim j$, and otherwise $i\nsim j$.  
%
%
%
%
Graphs in this paper are \emph{simple}, that is, they have at most one edge between a pair of nodes.  We consider two types of edges:  \emph{arrows} ($i\fra j$) and \emph{bidirected edges} or  \emph{arcs} ($i\arc j$). We do not consider graphs that have simultaneous third type of edge: \emph{undirected edges} or \emph{lines} ($i\ful j$).  The completely undirected graph, obtained from a graph $G$,  is called the \emph{skeleton} of $G$, denoted by $\sk(G)$; it is obtained by removing all arrowheads from the edges of the graph $G$. We say that we \emph{direct} the edges of a skeleton by assigning  arrowheads at the edges, which results in a graph.

A \emph{subgraph} of a graph $G_1$ is graph $G_2$ such that $V(G_2)\subseteq V(G_1)$ and $E(G_2)\subseteq E(G_1)$ and the assignment of endpoints to edges in $G_2$ is the same as in $G_1$. An \emph{induced subgraph} by nodes $A\subseteq V$ is a subgraph that contains all and only nodes in $A$ and all edges between two nodes in $A$.

A \emph{walk} is a list $\langle v_0,e_1,v_1,\dots,e_k,v_k\rangle$ of nodes and edges such that for $1\leq i\leq k$, the edge $e_i$ has endpoints $v_{i-1}$ and $v_i$. A \emph{path} is a walk with no repeated node or edge. When we define a path, we only write the nodes (and not the edges). 
A \emph{cycle} is a walk  with no repeated nodes or edges except for $v_0=v_k$.

We call the first and the last nodes \emph{endpoints} of the path and all other nodes \emph{inner nodes}. A path can also be seen as a certain type of connected subgraph of $G$;  a \emph{subpath} of a path $\pi$ is  an induced connected subgraph of $\pi$.

For an arrow $j\fra i$, we say that the arrow is \emph{from} $j$ \emph{to} $i$. We also call $j$ a \emph{parent} of $i$, $i$ a \emph{child} of $j$ and we use the notation $\pa(i)$ for the set of all parents of $i$ in the graph.
In the cases of $i\fra j$ or $i\arc j$ we say that there is an arrowhead at $j$ or pointing to $j$. 
A path $\langle i=i_0,i_1,\dots,i_n=j\rangle$ is \emph{directed} from $i$ to $j$ if all $i_ki_{k+1}$ edges are arrows pointing from $i_k$ to $i_{k+1}$. If there is a directed path from $j$ to $i$, then node $j$ is an \emph{ancestor} of $i$ and $i$ is a \emph{descendant} of $j$. We denote the set of ancestors of $i$ by $\an(i)$;  unlike some authors, we do not assume that $i\in\an(i)$.

A \emph{tripath} is a path with three nodes. Note that \citet{sad13} used the term V-configuration for such a path. We follow \citet{kii84} and most texts by letting a \emph{V-configuration} be a tripath with non-adjacent endpoints.
The inner node $t$ in each of the  three tripaths
$$i\fra\,
t\fla\,j, \  i\arc\,t\fla\, \ j, \  i\arc\,t\arc\,j$$
is
a \emph{collider} (or a collider node) and the inner node of any other tripath
is a \emph{non-collider} (or a non-collider node) on the tripath or, more generally, on any path of which the tripath is a subpath; i.e., a node is a collider if two arrowheads meet. Notice that in some causal inference literature, the term \emph{V-structure} is used for a collider V-configuration.  A path is called a \emph{collider} path if all its inner nodes are colliders.

In this paper, we are mainly concerned with  \emph{directed ancestral graphs};  these are \emph{ancestral graphs} (AGs) \citep{ric02} without lines, i.e., graphs with arrows and arcs without a directed cycle and without an arc with one endpoint being an ancestor of the other.
Sometimes we  simply refer to such graphs as ``graphs," dropping any modifiers.   Some results we use in this paper have originally been proven for the more general classes of chain mixed graphs (CMGs) and anterial graphs (AnGs) \citep{sad16} as well as acyclic directed mixed graphs (ADMGs) \citep{ric03}, all of which contain directed ancestral graphs.

The class of directed ancestral graphs also contains 
\emph{bidirected graphs} (BGs), containing  only bidirected edges, and \emph{directed acyclic graphs} (DAGs), containing only arrows and being acyclic. DAGs have been particularly useful  to describe causal Markov relations; see for example \citet{kii84,pea88,lau88,gei90,spioo}.
For an extensive discussion on the subclasses of acyclic graphs and their relationships and hierarchy, see \citet{sadl16}.

\subsection{Independence models and their properties}\label{sec:propind}
An \emph{independence model} $\mathcal{J}$, over a finite set $V$, is a set of triples $\langle X,Y\cd Z\rangle$, called \emph{independence statements}, where $X$, $Y$, and $Z$ are disjoint subsets of $V$; $Z$
may be empty, but $\langle \varnothing,Y\cd Z\rangle$ and $\langle X,\varnothing\cd Z\rangle$ are always included in $\mathcal{J}$. The independence statement $\langle X,Y\cd Z\rangle$ is read as ``$X$ is independent of $Y$ given $Z$.'' Independence models may  have a  probabilistic interpretation, but this is not necessarily the case. Similarly, not all independence models can be easily represented by graphs. For further discussion on general independence models, see \citet{stu05}.

In order to define probabilistic independence models, consider a set $V$ and a collection of random variables
$\{X_\alpha\}_{\alpha\in V}$ with state spaces $\{\mathcal{X}_\alpha\}_ {\alpha\in V}$ and joint distribution $P$. We let $X_A=\{X_v\}_{v\in A}$  for each subset $A$ of $V$. For disjoint subsets $A$, $B$, and $C$ of $V$
we use the short notation $A\cip B\cd C$ to denote that $X_A$ is \emph{conditionally independent of $X_B$ given $X_C$} \citep{daw79,lau96}, that is,  for any measurable $S\subseteq \mathcal{X}_A$ and $P$-almost all $x_B$ and $x_C$, we have
$$\mathbb{P}(X_A \in S\cd X_B=x_B, X_C=x_C)=\mathbb{P}(X_A \in S\cd X_C=x_C),$$
where we make suitable interpretations in terms of regular conditional probabilities when necessary; see for example \citet{MR1484954}.
We can now induce an independence model $\mathcal{J}(P)$ by letting
\begin{displaymath}
\langle A,B\cd C\rangle\in \mathcal{J}(P) \iff  A\ci B\cd C \text{ w.r.t.\ $P$}.
\end{displaymath}
Similarly we use the notation $A\notci B\cd C$ for $\langle A,B\cd C\rangle\notin \mathcal{J}(P)$. In this paper, we present the results for ``$P$," but in fact all the conditions and assumptions are on $\mathcal{J}(P)$. We remark that all the results and conditions could be presented using an abstract independence model $\mathcal{J}$.

In order to use graphs to represent independence models, the notion of \emph{separation} in a graph is fundamental.    A path $\pi$ is \emph{connecting}  given $C$ if every  collider node of $\pi$ is in $C$ or an ancestor of a node in $C$, and every non-collider node is not
in
$C$.    For pairwise disjoint subsets $(A,B,C)$, we write  $A\dse B\cd C$ if there are no connecting paths between $A$ and $B$ given $C$, and say that  $A$ and $B$ are \emph{separated given} $C$.  We will also use the notation  $A\notdse B\cd C$ for $A$ and $B$ not separated given $C$.
%
%
%
%
 This separation criterion is generally known as the \emph{$m$-separation}, used in \citet{ric02}, \citet{wer11}, and \citet{wers11}, for ancestral graphs and summary graphs, and is a generalization of the  \emph{d-separation} of \citet{pea88}. 

If $A$, $B$, or $C$ has only one member $\{i\}$, $\{j\}$, or $\{k\}$, for
   better readability, we write $\langle i,j\cd k\rangle\in \mathcal{J}$ instead of $\langle \{i\},\{j\}\cd \{k\}\rangle\in \mathcal{J}$; and similarly for $i\dse j\cd k$  and $i\ci j\cd k$. We also write $A\dse B$ when $C=\varnothing$; and similarly $A\ci B$.

A graph $G$ induces an independence model $\mathcal{J}(G)$ by separation, letting
\begin{displaymath}
\langle A,B\cd C\rangle\in \mathcal{J}(G)\iff A\dse B\cd C \text{ in } G.
\end{displaymath}

A graph is called \emph{maximal} if the absence  of an edge between $i$ and $j$ corresponds to a conditional separation statement for $i$ and $j$, i.e., there exists for some $C$ a statement of form $i\dse j\cd C$. In addition, we call two graphs $G$ and $H$
 \emph{Markov equivalent} if $\mathcal{J}(G)=\mathcal{J}(H)$. Conditions for Markov equivalence for most classes of graphs are known; see \citet{ver90,ali09,wers11}. Notice that two Markov equivalent maximal graphs have the same skeleton. We require conditions for Markov equivalence of the classes of graphs discussed here.    A collider path $\pi=\langle i,B,j\rangle$ is called a \emph{minimal collider path} if $i\nsim j$ and there is no proper subset $B'\subset B$ such that $\langle i,B',j\rangle$ is a collider path between $i$ and $j$.
\begin{lemma}[\cite{zha04}]
	\label{lem:Mark-eq-MAG}
Two maximal ancestral graphs (and consequently maximal directed ancestral graphs) are Markov equivalent if and only if they have the same skeletons and minimal collider paths.
\end{lemma}
For DAGs, there is a simpler condition for Markov equivalence:
\begin{lemma}[\cite{ver90}]
	\label{lem:Mark-eq-DAG}
Two DAGs are Markov equivalent if and only if they have the same skeletons and collider V-configurations.
\end{lemma}

An independence model $\mathcal{J}$ over a set $V$ is a \emph{semi-graphoid} if it satisfies the four following properties for disjoint subsets $A$, $B$, $C$, and $D$ of $V$:
 \begin{enumerate}
    \item $\langle A,B\cd C\rangle\in \mathcal{J}$ if and only if $\langle B,A\cd C\rangle\in \mathcal{J}$ (\emph{symmetry});
    \item if $\langle A,B\cup D\cd C\rangle\in \mathcal{J}$, then $\langle A,B\cd C\rangle\in \mathcal{J}$ and $\langle A,D\cd C\rangle\in \mathcal{J}$ (\emph{decomposition});
    \item if $\langle A,B\cup D\cd C\rangle\in \mathcal{J}$, then $\langle A,B\cd C\cup D\rangle\in \mathcal{J}$ and $\langle A,D\cd C\cup B\rangle\in \mathcal{J}$ (\emph{weak union});
    \item if $\langle A,B\cd C\cup D\rangle\in \mathcal{J}$ and $\langle A,D\cd C\rangle\in \mathcal{J}$,
    then $\langle A,B\cup D\cd C\rangle\in \mathcal{J}$ (\emph{contraction}).
 \end{enumerate}
Notice that the reverse implication of contraction clearly holds by decomposition and weak union. A semi-graphoid for which the reverse implication of the weak union property holds is said to be a \emph{graphoid}; that is, it also satisfies
\begin{itemize}
	\item[5.] if $\langle A,B\cd C\cup D\rangle\in \mathcal{J}$ and $\langle A,D\cd C\cup B\rangle\in \mathcal{J}$, then $\langle A,B\cup D\cd C\rangle\in \mathcal{J}$ (\emph{intersection}).
\end{itemize}
Furthermore, a graphoid or semi-graphoid for which the reverse implication of the decomposition property holds is said to be \emph{compositional}, that is, it also satisfies
\begin{itemize}
	\item[6.] if $\langle A,B\cd C\rangle\in \mathcal{J}$ and $\langle A,D\cd C\rangle\in \mathcal{J}$, then $\langle A,B\cup D\cd C\rangle\in \mathcal{J}$ (\emph{composition}).
\end{itemize}
Another important property of independence models is \emph{singleton-transitivity} (also called \emph{weak transitivity} in \citet{pea88}). For $i$, $j$, and $k$, single elements in $V$,
 \begin{itemize}
	\item[7.] if $\langle i,j\cd C\rangle\in \mathcal{J}$ and $\langle i,j\cd C\cup\{k\}\rangle\in \mathcal{J}$, then $\langle i,k\cd C\rangle\in \mathcal{J}$ or $\langle j,k\cd C\rangle\in \mathcal{J}$ (singleton-transitivity).
\end{itemize}

Probabilistic independence models are always semi-graphoids \citep{pea88}, whereas the converse is not necessarily true; see \citet{stu89}. If $P$ has strictly positive density, the induced independence model is always a graphoid; see, for example, Proposition 3.1 in \citet{lau96}. See also \citet{pet15} for a necessary and sufficient condition for the intersection property to hold. If the distribution $P$ is a regular multivariate Gaussian distribution, then  $\mathcal{J}(P)$ is a singleton-transitive compositional graphoid; for example see \citet{stu05} and \cite{pea88}. For this reason, a different axiomatization of singleton-transitive compositional graphoid is called a \emph{Gaussoid} in \citet{lne07}. On the other hand, separation in graphs satisfies all these properties; see Propositions 1 and 2 in \citet{sad17}.
If $\mathcal{J}(P)$ satisfies any of these properties, then we may simply say that ``P satisfies that property."

\subsection{Ordered upward- and downward-stabilities}\label{sec:ouds} The definitions and results below are a special case of those in Section 4 of \citet{sad17}.   Consider a \emph{partial order} $\leq$ over the set $V$. If $a \leq b$ or $b \leq a$, then $a$ and $b$ are \emph{comparable}; otherwise they are \emph{incomparable}.  Henceforth, we speak of ``orders" referring to  ``partial orders." 

We say that a  graph $G=(V,E)$ admits a \emph{valid order} $\leq$ if for nodes $i$ and $j$ of $G$, $i\fra j$ implies that $i >j$; and $i\arc j$ implies that $i$ and $j$ are incomparable. Notice that this specifies the partial order via its cover relations.
There may be many different orders that are valid for a graph. However, we obtain a unique ordering for a graph as follows:   Given a graph $G$, if $i\notin\an(j)$ and $j\notin\an(i)$,
 then set $i$ and $j$ to be incomparable. Otherwise, let $i>j$ when $i\fra j$. We call this ordering the \emph{minimal order} for $G$ since it gives the fewest possible comparable pairs of nodes.
It can also be seen that if $G$ is ancestral, then the minimal order for $G$ is a valid order.

We now exploit the ordering for independence models in order to define two other properties in addition to the seven properties defined in Section \ref{sec:propind} (namely singleton-transitive compositional graphoid axioms). We say that an independence model $\mathcal{J}$ over the set $V$ respectively satisfies \emph{ordered upward-} and \emph{downward-stability} w.r.t.\ an order $\leq$ of $V$ if the following hold:
  \begin{itemize}
	\item[8.] if $\langle i,j\cd C\rangle\in \mathcal{J}$,  then $\langle i,j\cd C\cup\{k\}\rangle\in \mathcal{J}$ for every $k\in V\setminus\{i,j\}$ such that $i< k$ or $j< k$ (ordered upward-stability);
	\item[9.] if $\langle i,j\cd C\rangle\in \mathcal{J}$, then $\langle i,j\cd C\setminus\{k\}\rangle\in \mathcal{J}$ for every $k\in V\setminus\{i,j\}$ such that $i\nless k$, $j\nless k$, and $l\nless k$ for every $l\in C\setminus\{k\}$ (ordered downward-stability).
\end{itemize}
Sometimes we refer to the pair of Conditions 8 and 9, together, simply as \emph{ordered stability}.

Ordered upward-stability is a generalization of a modification of \emph{upward stability}, defined in \citet{fal15}, and \emph{strong union}, defined in \citet{pea85}, for undirected graphs. A similar concept was used in \citet{cla12} for the class of DAGs for inferring causal relations; see Lemma 2 of the mentioned manuscript.

\citet{sad17} proved  that the independence model induced by ancestral graphs satisfies ordered upward- and downward-stability w.r.t.\ their minimal order. If an independence model $\mathcal{J}$ satisfies ordered upward- and downward-stability w.r.t.\ the minimal order of a graph $G$,
then we simply say that ``$\mathcal{J}$ satisfies upward- and downward-stability w.r.t.\ $G$."  In addition, similar to compositional graphoids, if $\mathcal{J}(P)$ satisfies ordered upward- or downward-stability, then we may simply say that ``$P$ satisfies that property."

\subsection{Markov and faithful independence models}\label{sec:prob}
For a graph $G=(V,E)$, a probability distribution defined over $V$ satisfies the \emph{global Markov property} w.r.t.\ $G$, or is simply \emph{Markovian} to $G$, if for  disjoint subsets $A$, $B$, and $C$ of $V$ it holds that  $$A\dse B\cd C \implies A\ci B\cd C,$$
or equivalently $\mathcal{J}(G)\subseteq\mathcal{J}(P)$.
Notice that every probability distribution over $V$ is Markovian to the complete graph with the node set $V$.

We say that $P$ and a graph $G$ are
\emph{faithful}
if $\mathcal{J}(P)= \mathcal{J}(G)$. If $P$ and $G$ are faithful, then we may sometimes also say that $P$ is \emph{faithful to}  $G$ or vice versa, although in principle faithfulness is a symmetric relation.
Thus, if $P$ and $G$ are faithful, then  in addition to $P$ being Markovian to $G$, every conditional independence statement corresponds to a separation in $G$. Notice that, originally in \citet{spioo}, faithfulness was defined only to be the opposite direction to the Markov property.  
If there is a graph $G$ and a distribution $P$ that are faithful, then we say that $P$ is \emph{graphical} and $G$ is \emph{probabilistic}.

For a given probability distribution $P$, we define the \emph{skeleton} of $P$, denoted by $\sk(P)$, to be the undirected graph that is obtained from $P$ as follows:  we define the node set of $\sk(P)$ to be $V$, and for every pair of nodes $i,j$, we check whether $i\ci j\cd C$  holds for some $C\subseteq V\setminus\{i,j\}$; if it does not then we draw an edge between $i$ and $j$. The resulting undirected graph is the same as the one obtained by the first step of the SGS algorithm \citep{gly87}. 

Suppose that there exists an order $\leq$ over $V$. We can uniquely direct the edges of $\sk(P)$ based on $\leq$ in order to define a $G(P,\leq)$ induced by $P$ and $\leq$. It can be seen that $G(P,\leq)$ is ancestral.

We are only interested in orderings that are minimal orders of graphs. Given $P$, we define an ordering $\leq$ to be \emph{$P$-compatible} if $\leq$ is the minimal order for $G(P,\leq)$.
%
%
%

If $P$ is Markovian to a graph $G$, then $\sk(P)$ is a subgraph of $\sk(G)$; see \citet{sad17}. Hence, if $P$ is Markovian to a graph $G$ such that $\sk(G)=\sk(P)$, then $G$ has the fewest number of edges among those to which $P$ is Markovian. We say that $P$ is \emph{minimally Markovian} to a graph $G$  if $P$ is Markovian to $G$ and $\sk(G)=\sk(P)$ \citep{sad17}.  Notice that only minimally Markovian independence models to a graph can also be faithful to the graph.
\begin{remark}[Related versions of  minimality]
\label{rem:minimal}
There are several
similarly defined
minimality assumptions in the literature.
Notice that we do not
define the concept of minimal Markovian as an assumption in this paper, since minimality is implied by some other assumptions used in this paper; see Corollary \ref{coro:minmark}.

The definition given here is subtly different from the
{(causal) minimality assumption}
\citep{pea09,spioo,zha08,nea04}, which holds if $P$ is Markovian to $G$, but not to any proper subgraph of $G$ --- in fact, it can be seen that being minimal Markovian implies causal minimality.  It is also different from the
{P-minimality}
assumption \citep{zha13}, which states that no proper
{I-structure}
of $G$ is Markovian to $P$, where I-structure of $G$ is any subgraph of $G$ whose induced independence model is a subset of that of $G$. Again,  being minimal Markovian implies P-minimality since there could be a graph with fewer edges than $G$ that is not an $I$-structure of $G$.

The
{sparsest Markov representation}
(SMR) assumption \citep{ras18} for $G$ and $P$ (originally defined for DAGs) states that  $P$ is Markovian to $G$ and  every $G'$ to which $P$ is Markovian and not Markov equivalent to $G$ has more edges than $G$ does.  A similar concept is the
{frugality}
assumption \citep{for18}. In the special case of interest here, where $G$ is Markovian to $P$,  it states that $G$ satisfies frugality if it is in the set of sparsest Markov representations.  If $G$ and $P$ satisfy SMR,
 then $G$ is not necessarily minimally Markovian to $P$ (see Figure 5 of \citet{for18} as an example).
It is also possible to create two non-Markov equivalent graphs with the same skeleton, and a $P$ that is Markovian to both. This implies that the minimal Markovian assumption does not imply SMR either.
\end{remark}
We make use of the necessary and sufficient conditions for faithfulness to directed ancestral graphs,  proven in \citet{sad17} for the more general case of anterial graphs:
\begin{proposition}[\cite{sad17}]
	\label{thm:2}
Let $P$ be a probability distribution defined over $X_V$.  Then  $P$ is graphical if and only if
\begin{enumerate}[(A)]
  \item
\label{sing}
$P$ is a singleton-transitive compositional graphoid; and
  \item
  \label{down}
   there exists a $P$-compatible order $\lesssim$ over $V$ w.r.t.\ which $P$ satisfies ordered downward- and upward-stability.
\end{enumerate}
In addition, if \eqref{sing} and \eqref{down} are satisfied, then $P$ is faithful to $G(P,\lesssim)$.
\end{proposition}

\subsection{Structural causal models}\label{sec:SCM}
Here we define the \emph{structural causal models} (SCMs) (also known as the \emph{structural equation models} (SEMs)) \citep{pea09,spioo} for the class of directed ancestral graphs, as introduced in \citet{ric02}. Consider a graph $G$ with $N$ nodes, which in the context of causal inference may be referred to as the ``true causal graph." An SCM $\mathfrak{C}$  associated with $G$ is defined as a collection of $N$ equations
\begin{equation}\label{eq:scm}
X_i=\phi_i(X_{\pa(i)},\epsilon_i), \nn i\in\{1,\dots, N\},
\end{equation}
where $\pa(i)$ is defined on $G$, and, for any subsets $A$ and $B$, we require that $\epsilon_A \ci \epsilon_B$  if and only if, in $G$, there is no arc between any node in $A$ and any node in $B$. Notice that in the more widely-used case where $G$ is a DAG, all
the
$\epsilon_i$ are jointly independent. Sometimes  the $X_i$ are called \emph{endogenous variables} and the \emph{noises} $\epsilon_i$  are also called \emph{exogenous variables} (which are generally assumed to be latent). For both mathematical and causal discussions on SCMs with DAGs, see \citet{pet17}. In the case where
$\phi_i$
represent multivariate linear regressions, SCMs are also known as the system of linear equations \citep{cox96}.

Notice that the minimal order of $G$ induces an order for $X_i$s in $\mathfrak{C}$. In addition, by using the equations iteratively, it is easy to see that each $X_i$ can be written as a function of $\epsilon_{\an(i)\cup\{i\}}$. An SCM also entails a unique joint distribution $P$
over the variables $(X_1,\dots,X_N)$.

A more complete representation of the graph associated with an SCM is obtained by including the random variables $\epsilon_i$. These graphs, called \emph{augmented graphs},
are
denoted by $G(X,\epsilon)$. We construct the augmented graph in our setting as follows: we denote the nodes in $G$ by the original notation $\{X_1,\dots,X_N\}$; we add the nodes $\{\epsilon_1,\dots,\epsilon_N\}$ and add arrows from each $\epsilon_i$ to $X_i$; we also remove all arcs between $X_i$s, and instead of every $X_iX_j$-arc, add an arc between $\epsilon_i$ and $\epsilon_j$. 

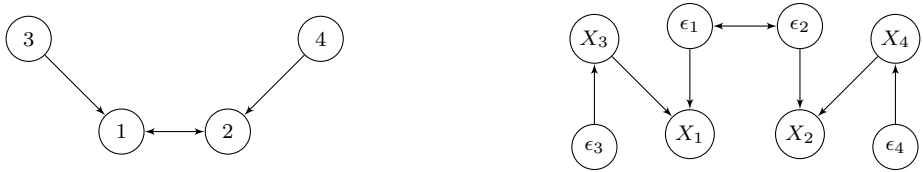
\begin{figure}[htb]
\centering
\begin{tikzpicture}[node distance = 8mm and 8mm, minimum width = 6mm]
    \begin{scope}
      \tikzstyle{every node} = [shape = circle,
      font = \scriptsize,
      minimum height = 6mm,
      inner sep = 2pt,
      draw = black,
      fill = white,
      anchor = center],
      text centered]
      \node(1) at (0,0) {$3$};
      \node(2) [ below right = of 1] {$1$};
      \node(3) [right = of 2] {$2$};
      \node(4) [above right = of 3] {$4$};
      \node(1n) [right = 30mm of 4] {$X_3$};
      \node(2n) [ below right = of 1n] {$X_1$};
      \node(3n) [right = of 2n] {$X_2$};
      \node(4n) [above right = of 3n] {$X_4$};
      \node(5n) [below = of 1n] {$\epsilon_3$};
      \node(6n) [above = of 2n] {$\epsilon_1$};
      \node(7n) [above = of 3n] {$\epsilon_2$};
      \node(8n) [below = of 4n] {$\epsilon_4$};
    \end{scope}
		
    \begin{scope}
    \tikzstyle{every node} = [node distance = 6mm and 6mm, minimum width = 6mm,
    font= \scriptsize,
      anchor = center,
      text centered]

\end{scope}

    \begin{scope}[->, > = latex']
    \draw (1) -- (2);
    \draw (4) -- (3);
	\draw (5n) -- (1n);
	\draw (1n) -- (2n);
    \draw (8n) -- (4n);
    \draw (4n) -- (3n);
    \draw (6n) -- (2n);
    \draw (7n) -- (3n);
    \end{scope}
    \begin{scope}[<->, > = latex']
    \draw (2) -- (3);
    \draw (6n) -- (7n);
    \end{scope}

    \end{tikzpicture}
		\caption{{\footnotesize A graph $G$ (left) and  its augmented graph $G(X,\epsilon)$ (right).}}\label{fig:exaug}
		\end{figure}

For example, in Figure \ref{fig:exaug}, the graph on the right, is the augmented graph of the, graph on the left. An SCM associated to these graphs, in which $\phi_i$ are linear is as follows:
$$X_4=\epsilon_4,\nn X_3=\epsilon_3,\nn X_2=\alpha X_4+\epsilon_2, \nn X_1=\beta X_3+\epsilon_1.$$

The above setting and construction of the augmented graph is in line with the idea originally presented in \citet{ric02} for linear functions and Gaussian distributions, although unlike in the mentioned manuscript, noises corresponding to endpoints of variables adjacent by an arc are \emph{required} to  be dependent. This setting is different from the dominant setting in the literature (e.g.\ \citet{pet17}) for directed ancestral graphs, where for every arc $X_iX_j$ it is assumed that there exists an exogenous variable $\epsilon_{ij}$, where, in the directed acyclic augmented graph, it points to both $X_i$ and $X_j$, and that $\epsilon_{ij}$s are mutually independent.

\begin{remark}\label{rem:joint-noise} Our setting is more general than the other alternative presented
above.
Notice that we can write a noise $\epsilon_i$ as a tuple $(\epsilon_{is})_{s\in S}$, where $S$ is the set of nodes connected to $i$ by an arc. However, by writing $\epsilon_{ij}=(\epsilon_i,\epsilon_j)$, we would not obtain independent noises in the other setting.
\end{remark}
\section{General setting for structure learning}\label{sec:gen-set}
We start by assuming the existence of an \emph{unknown} distribution $P$, which is Markovian to an unknown \emph{(true) causal graph} $G_0$. Indeed, this setting requires additional models and assumptions to be considered ``causal."
This stems from the fact that, in the causal model, interventional probabilities and observational conditional probabilities coincide.
For more discussion on this, see \citet{pea09,spioo,pet17}. Nevertheless, the setting presented here is what is required for causal discovery using conditional independence under any causal model.

The main goal of structure learning is to use observational data to find a graph that belongs to the Markov equivalence class of $G_0$. The distribution $P$ induces an independence model $\mathcal{J}(P)$ whose elements we can test for using the observed data. Although we can, in principle, test each statement in $\mathcal{J}(P)$, testing for all independence statements is computationally intractable. In a constraint-based structural learning algorithm, one uses some of the independence statements in  $\mathcal{J}(P)$ to construct a graph, which is (hopefully) in the Markov equivalence class of $G_0$. Notice that, in this paper, we concern ourselves only with the oracle setting, i.e., assuming we have an oracle that will always correctly tell us whether or not a given conditional independence relation holds in the distribution. 

Indeed the exact relationship between $P$ (or more specifically $\mathcal{J}(P)$) and the constructed graph depends on the specific algorithm used.  Here, we propose the following class of generic algorithms.  

\begin{definition}[Natural structure-learning algorithm]
	\label{algorithm}  Let $P$ be a joint probability distribution on $X_V$.
We say that a natural constraint-based structure-learning algorithm is an algorithm whose output graph with node set $V$, denoted by $G(P)$, has the following properties:
\begin{enumerate}
  \item It holds that $\sk(G(P))=\sk(P)$;
  \item The  distribution $P$ satisfies ordered downward- and upward-stability w.r.t.\ $G(P)$, that is, $G(P)=G(\mathcal{J}(P),\leq)$, for a partial order $\leq$, w.r.t.\ which  $P$ satisfies ordered downward- and upward-stability.
\end{enumerate}
\end{definition}
In short, the diagram below shows the relationship between the concepts used in this setting:
$$G_0 \stackrel{\text{Markovian}}{\dashleftarrow} P \stackrel{\text{conditional independence}}{\rightarrow} \mathcal{J}(P)\stackrel{\text{Algorithm}}{\rightarrow} G(P) \stackrel{\text{Markov equivalent}}{\sim} G_0$$

In the same fashion as many constraint-based structural learning algorithms, the skeleton of the output of a natural structure-learning algorithm is the same as $\sk(P)$. Under the faithfulness assumption of $P$ and $G_0$, we will show in Section \ref{sec:faith} that not only is $G(P)$  Markov equivalent to $G_0$, but also $P$ satisfies ordered downward- and upward-stability w.r.t.\ every graph Markov equivalent to $G_0$; thus the class of natural structure-learning algorithms contains all well-known ``exact" algorithms that work under faithfulness, so we use  term ``natural.''    Under alternative conditions, provided in Section \ref{sec:alg}, we show that the output $G(P)$ is Markov equivalent to $G_0$.   Notice again that we will always use $G(P)$ to denote an output of a natural structure-learning algorithm.

The algorithms in the literature, generally generate a completed partially directed acyclic graph (CPDAGs) \citep{and97,spioo} or  partial ancestral graphs (PAGs) \citep{ric02,ali09}, which are a representative of the Markov equivalence class of $G_0$ with the fewest arrowheads. Here, we are not concerned with making a  CPDAG or a PAG from $G(P)$.

We assume, without loss of generality, that the causal graph $G_0$ is maximal: the reason is that  if $G_0$ is not maximal, then there is always a maximal graph that is Markov equivalent to $G_0$, which could be used to find graphs Markov equivalent to $G_0$; see \citet{ric02,sadl14}.  Note that the first condition of the natural structure-learning algorithm in Definition \ref{algorithm} also ensures that the generated $G(P)$ is maximal.

Since Markov equivalent maximal graphs have the same skeleton, the obtained $G(P)$ must have the same skeleton with $G_0$. We are looking for conditions, under which this holds.

Let $\an(i,j)=(\an(i)\cup\an(j))\setminus\{i,j\}$. We say that $P$ is \emph{converse pairwise Markovian} to $G_0$ if
\begin{equation}
	\label{converse-def}
i\ci j\cd \an(i,j)\Rightarrow i\nsim j \text{ in } G_0;
\end{equation}
note that the converse of this condition is the well-known pairwise Markov property \citep{lau96}. To ensure the equality of the mentioned skeletons, we need the following lemma:
\begin{lemma}\label{lem:skel}
Suppose that the distribution $P$ is Markovian to
the causal graph
 $G_0$, and that the following two conditions are satisfied:
\begin{enumerate}[(a)]
  \item
  \label{condition1}
  $P$ is converse pairwise Markovian to $G_0$;
  \item
  \label{condition2}
  if there exists $C$, not containing $i$ and $j$,  such that  $i\ci j\cd C$, then $i\ci j\cd \an(i,j)$.
\end{enumerate}
Then $\sk(G_0)=\sk(G(P))$, where %
 $G(P)$ is an  output of a natural structure-learning algorithm.
\end{lemma}
\begin{proof}
	We have on $G_0$ the equivalence
$$i\nsim j \text{ in } G_0 \iff  i\ci j\cd \an(i,j),$$	
	where assumption \eqref{condition1} gives the reverse direction, and the forward direction follows from the fact that  $P$ is Markovian to $G_0$, and $\an(i,j)$ always separates non-adjacent $i$ and $j$.  Since $\sk(P)=\sk(G(P))$, we have  the equivalence
$$i\nsim j \text{ in } G(P)\iff \exists C \text{ s.t.\ } i\ci j\cd C.$$
Hence, assumption \eqref{condition2} allows us to tie these two conditions together and obtain that
$$i\nsim j \text{ in } G(P)  \iff i\nsim j \text{ in } G_0,$$
from which the desired  result follows.
\end{proof}
\begin{remark}
The two conditions of Lemma \ref{lem:skel} are equivalent to the
{adjacency-faithfulness}
condition \citep{ram06,zha08} (defined for DAGs), which states that for every edge between $i$ and $j$ in $G_0$, there are no independence statements $i\ci j \cd C$ for any $C$. This breakdown of the adjacency faithfulness assumption allows us to later provide useful sufficient conditions (in general and for SCMs) that ensure that these two conditions are satisfied; see also the discussion below.
\end{remark}
We need some additional assumptions on the distributions to determine when the assumptions of Lemma \ref{lem:skel} hold; i.e., when an existing edge implies dependence given the joint ancestors, and when an arbitrary conditional independence of a pair implies their conditional independence given their joint ancestors. A condition under which these assumptions hold is faithfulness, as will be discussed in the next section. However, weaker conditions will be provided in Section \ref{sec:alg}. Specific conditions for SCMs will also be provided in Section \ref{sec:st-l-SCM}.
\section{Structure learning under faithfulness}
\label{sec:faith}
Faithfulness of $P$ and $G_0$ is a main assumption in the literature for structure learning. It is well-known that under faithfulness, $P$ is \emph{identifiable}, i.e., an ancestral graph (or a DAG) (in reality a PAG (or a CPDAG)) Markov equivalent to $G_0$ can be identified \citep{spioo,zha08}. Under faithfulness, if $P$  is faithful to $G_0$,  then any structure-learning algorithm should provide an
output $G(P)$ that is faithful to $P$, in the setting where  an oracle can be queried regarding the validity of statements of conditional independence.

By Proposition \ref{thm:2}, if the assumption of $P$ being graphical holds, then $P$ is singleton-transitive compositional graphoid and there is an ordering $\leq$ of the nodes of $\sk(P)$ w.r.t.\ which $P$ satisfies ordered downward- and upward-stability. Notice that among all these conditions required for graphicality, ordered downward- and upward-stability are the only ones that are related to, and in fact determine, the direction of the edges of $\sk(P)$; these considerations also motivate the definition of natural structure-learning algorithms in Definition \ref{algorithm}.

We work towards a proof, that under faithfulness, a natural  structure-learning algorithm works.


\begin{proposition}\label{prop:skel}
Suppose that $P$ is faithful to
the causal graph
 $G_0$.
The output  of a natural structure-learning algorithm satisfies $\sk(G_0)=\sk(G(P))$.
\end{proposition}
{\bf Proof\ }
We show that the conditions of Lemma \ref{lem:skel} hold.
 \begin{enumerate}[(a)]
 	\item
 If $i\ci j\cd \an(i,j)$,  then  faithfulness gives, $i\dse j\cd \an(i,j)$ in $G_0$, which implies that $i\nsim j$ in $G_0$.
\item
 If there exists a $C$ such that $i\ci j\cd C$,  then   faithfulness gives,  $i\dse j\cd C$ in $G_0$, which  implies that $i\nsim j$ in $G_0$. Since $P$ is Markovian to $G_0$, we have that $i\ci j\cd \an(i,j)$.  \hfill \BlackBox
 \end{enumerate}

\begin{proposition}\label{thm:20}
Suppose that $P$ is faithful to
the causal graph
$G_0$. All potential outputs of natural structure-learning algorithms are Markov equivalent, and they are Markov equivalent to $G_0$. 
\end{proposition}
\begin{proof}
Since $P$ is graphical, by Proposition \ref{thm:2}, it is singleton-transitive compositional graphoid. Let $G(P)$ and $H(P)$ be two potential graphs generated by a natural structure-learning algorithm. By Proposition \ref{prop:skel}, we have  that $\sk(G_0)=\sk(G(P))=\sk(H(P))$.

In addition, $P$ satisfies ordered upward- and downward-stability w.r.t.\ both $G(P)$ and $H(P)$. Hence, again by Proposition \ref{thm:2}, $P$ is faithful to both. Therefore, these graphs are all Markov equivalent.
\end{proof}

In addition, the other direction of the above result also holds:
\begin{proposition}\label{thm:200}
Suppose that $P$ is faithful to $G_0$. If $G$  is Markov equivalent to $G_0$, then   $P$ satisfies ordered upward- and downward-stability w.r.t.\ $G$. 
\end{proposition}
\begin{proof}
The graph $G$ and $P$ have to be faithful. Therefore,  by Proposition \ref{thm:2}, we know that $P$ satisfies ordered upward- and downward-stability w.r.t.\ $G$.
\end{proof}

\begin{remark}
	\label{must}
Therefore, ordered downward- and upward-stability w.r.t.\ $G(P)$ are also necessary conditions of an output of a structure-learning algorithm that generates Markov equivalent graphs to $G_0$.
This justifies Definition \ref{algorithm}, when $P$ is graphical.
Hence, any ``exact"  (as opposed to approximate) constraint-based algorithms that work under the assumption of faithfulness (such as the SGS and PC algorithms \citep{spioo}) must (and, in fact, do) ensure that, w.r.t.\ their output, $P$ satisfies ordered downward- and upward-stability. 
\end{remark}

\section{Conditions and assumptions for the structure-learning algorithm}\label{sec:alg}


Recall that the general goal of constraint-based structure learning is to consider algorithms for which  the output $G(P)$ is Markov equivalent to (the unknown) causal graph $G_0$.
In Section \ref{sec:faith}, we saw that under the assumption of faithfulness of $P$ and $G_0$, this goal is achieved
for the class of algorithms falling under Definition \ref{algorithm},  the class of natural structure-learning algorithms.   However, in fact, faithfulness is not necessary for this goal, as will be shown below in Example \ref{ex:10}; see also Example \ref{ex:1} for when one focuses solely on the class of DAGs.
\begin{eexample}[Ordered stability and faithfulness]
	\label{ex:10}
Consider the graph in Figure \ref{fig:exuniquegen}, and assume that it is $G_0$ in the setting of this manuscript. Let $\mathcal{J}(P)$ be the set of statements implied by the global Markov property plus $i\ci j\cd k$. First of all, notice that $P$ does not satisfy singleton-transitivity. Notice also that $P$ is minimally Markovian to $G_0$, but $P$ and $G_0$ are not faithful. In addition, $P$ satisfies ordered upward- and downward-stability w.r.t.\ $G_0$.

However, all graphs, including the true causal graph $G_0$, w.r.t.\ which $P$ satisfies ordered downward- and upward-stability are Markov equivalent.  To see this,  suppose that there is a graph $H$ with the same skeleton, w.r.t.\ which $P$ satisfies ordered upward- and downward-stability, and that $H$ is not Markov equivalent to $G_0$. If, in $H$, $k$ is non-collider then $k>\ell$ or $k>i$, which, by ordered upward-stability, implies $i\ci \ell\cd k$; but this is not in $\mathcal{J}(P)$. Now, if in $H$, the node $\ell$ is collider, then $k,j>\ell$, which by ordered downward-stability and conditional independence $k\ci j\cd \ell$, implies $k\ci j$, which is again not  in $\mathcal{J}(P)$. 
\begin{figure}[htb]
\centering
\begin{tikzpicture}[node distance = 8mm and 8mm, minimum width = 6mm]
    \begin{scope}
      \tikzstyle{every node} = [shape = circle,
      font = \scriptsize,
      minimum height = 6mm,
      inner sep = 2pt,
      draw = black,
      fill = white,
      anchor = center],
      text centered]
      \node(1) at (0,0) {$i$};
      \node(2) [right = of 1] {$k$};
      \node(3) [right = of 2] {$\ell$};
      \node(4) [right = of 3] {$j$};
    \end{scope}
		
    \begin{scope}[->, > = latex']
    \draw (1) -- (2);
      \draw (3) -- (2);
    \draw (4) -- (3);
    \end{scope}

    \end{tikzpicture}
		\caption{{\footnotesize $G_0$}}\label{fig:exuniquegen}
		\end{figure}
\end{eexample}

Below, we find assumptions   under which  a natural structure-learning algorithm works; these assumptions can be thought of as a relaxation of faithfulness.
\subsection{Equivalence of the skeletons of $G_0$ and $G(P)$}
We first need to provide assumptions so that $G_0$ and $G(P)$ have the same skeletons. One of the main assumptions we employ is that
$P$ satisfies ordered downward- and upward-stability w.r.t.\ $G_0$.

\begin{theorem}\label{thm:ske}
Consider a distribution $P$ that is Markovian to
the causal graph
$G_0$ satisfying the following conditions:
\begin{enumerate}[(i)]
	\item
	\label{condition-i}
$P$ is converse pairwise Markovian to $G_0$;
\item
\label{condition-ii}
$P$ satisfies ordered upward- and downward-stability w.r.t.\ $G_0$.
\end{enumerate}
Then $\sk(G_0)=\sk(G(P))$.
\end{theorem}
\begin{proof}
It suffices to show that assumption \eqref{condition2} of  Lemma \ref{lem:skel} is satisfied. 
Suppose $i\ci j\cd C$.  Consider all the nodes in $C\setminus \an(i,j)$. Since $P$ satisfies ordered
downward-stability w.r.t.\ $G_0$, one can remove the node $k_0$ with a lowest order in this set to obtain $i\ci j\cd C\setminus\{k_0\}$. Inductively, we will obtain $i\ci j\cd C\cap \an(i,j)$. Now, by the use of ordered upward-stability, we add other nodes in $\an(i,j)$ to the conditioning set, and obtain $i\ci j\cd \an(i,j)$.
\end{proof}
\begin{corollary}\label{coro:minmark}
Consider a distribution $P$ that satisfies the conditions of Theorem \ref{thm:ske}.  Then  $P$ is minimally Markovian to $G_0$.
\end{corollary}
\begin{remark}
Notice that there exist cases where there is a unique Markov equivalent class of graphs including $G_0$ to which the distribution $P$ is Markovian, but the skeleton of $G_0$ is not the same as the skeleton of $P$; see \citet{for18}, where such examples are used as a motivation for the assumption of frugality; see Remark \ref{rem:minimal}. In such cases, conditions of Theorem \ref{thm:ske} (or equivalently adjacency-faithfulness) are not satisfied. This suggests that there could be extensions to the setting of this paper.
\end{remark}

\subsection{Markov equivalence of $G(P)$ and $G_0$}
In this section, we will describe under what additional testable conditions 
$G_0$ and $G(P)$ are Markov equivalent. We provide the conditions and results both for the case where $G_0$ is a DAG and the structure-learning algorithm only generates DAGs, and also for the general case where $G_0$ is a directed ancestral graph and the structure-learning algorithm generates such graphs.

We formalize the problem using the following concept: We say that a distribution $P$ satisfies the \emph{uniqueness property} (up to Markov equivalence) when the following holds: maximal graphs with skeleton $\sk(P)$ w.r.t.\ which $P$ satisfies ordered upward- and downward-stability are Markov equivalent. For the special case of DAGs, we say that a distribution $P$ satisfies the \emph{DAG-uniqueness property} if the graphs concerned in the definition of the uniqueness property are restricted to DAGs.

The distribution $P$ satisfying the uniqueness property is equivalent to the potential outputs of the structure-learning algorithm being unique up to Markov equivalence. This
is an obvious requirement of the structure-learning algorithm so that the algorithm cannot generate two graphs that are not Markov equivalent.

In order to ensure that the output of a natural structure-learning algorithm  is correct, we need the following important lemma:

\begin{lemma}\label{lem:essen}
Consider a distribution $P$ that is Markovian to $G_0$. Suppose that the conditions of Theorem \ref{thm:ske} are satisfied; that is,   $P$ is converse pairwise Markovian to $G_0$ and $P$ satisfies ordered upward- and downward-stability w.r.t.\ $G_0$.  %
Let $G(P)$ be the output of a natural structure-learning algorithm.
\begin{enumerate}[(a)]
  \item
   If $P$ satisfies the uniqueness property, then $G(P)$ is Markov equivalent to $G_0$.
  \item
   If $G_0$ is a DAG, then  $P$ satisfying the DAG-uniqueness property implies that a DAG output $G(P)$ is Markov equivalent to $G_0$.
\end{enumerate}
\end{lemma}
\begin{proof}
The proof follows from Theorem \ref{thm:ske}.
\end{proof}
In what follows,  we provide and study explicit conditions under which the uniqueness property is satisfied.

By Corollary \ref{coro:minmark}, we know that, under the assumptions for $P$, which we have so far utilized, $P$ is minimally Markovian to $G_0$.  In Example \ref{ex:1}, we observe that, if we replace the faithfulness assumption with the assumption of being minimal Markovian, then the uniqueness property no longer holds. 
\begin{eexample}\label{ex:1}
Let $\mathcal{J}(P)=\{1\ci 4\cd \{2,3\},\ 1\ci 4\cd 3,\ 2\ci 3\cd 4\}$. For example, there is an SCM that induces such an independence model: Let $\epsilon_1,\epsilon_2,\epsilon_3$, and $\epsilon_4$ be i.i.d. Define
$$X_1=\max(X_2,X_3,\epsilon_1);\nn X_2=\max(X_4,\epsilon_2);\nn X_3=\max(2X_4,\epsilon_3);\nn X_4=\epsilon_4.$$
It can be seen that $(X_i)_{i\in\{1,\dots,4\}}$ induces $\mathcal{J}(P)$.

First notice that $P$ does not satisfy the uniqueness property: w.r.t.\ both graphs $G_1$ and $G_2$ in Figure \ref{fig:exunique}(b) and \ref{fig:exunique}(c), $\mathcal{J}(P)$ satisfies ordered upward- and downward-stability. However, these two graphs are not Markov equivalent.


Since $P$ satisfies the compositional graphoid axioms (but not singleton-transitivity), this also shows that compositional graphoid plus ordered upward- and downward-stability are not in general sufficient for the uniqueness property to hold.

However, if we assume that the structure-learning algorithm is only DAG-generating, then there is a unique graph w.r.t.\ which  $P$ satisfies ordered upward- and downward-stability (and consequently the DAG-uniqueness property holds).
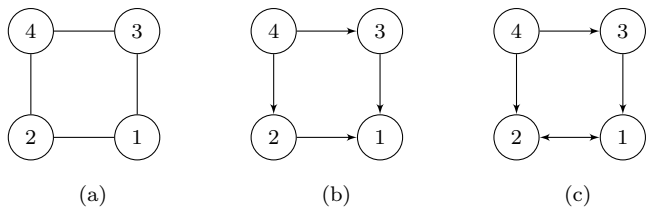
\begin{figure}[htb]
\centering
\begin{tikzpicture}[node distance = 8mm and 8mm, minimum width = 6mm]
    \begin{scope}
      \tikzstyle{every node} = [shape = circle,
      font = \scriptsize,
      minimum height = 6mm,
      inner sep = 2pt,
      draw = black,
      fill = white,
      anchor = center],
      text centered]
      \node(1) at (0,0) {$4$};
      \node(2) [right = of 1] {$3$};
      \node(3) [below = of 1] {$2$};
      \node(4) [below = of 2] {$1$};
      \node(1n) [right = 12mm of 2] {$4$};
      \node(2n) [right = of 1n] {$3$};
      \node(3n) [below = of 1n] {$2$};
      \node(4n) [below = of 2n] {$1$};
      \node(1nn) [right = 12mm of 2n] {$4$};
      \node(2nn) [right = of 1nn] {$3$};
      \node(3nn) [below = of 1nn] {$2$};
      \node(4nn) [below = of 2nn] {$1$};
    \end{scope}
		
    \begin{scope}
    \tikzstyle{every node} = [node distance = 6mm and 6mm, minimum width = 6mm,
    font= \scriptsize,
      anchor = center,
      text centered]
\node(a) [below right= 4mm of 3]{(a)};
\node(b) [below right= 4mm  of 3n]{(b)};
\node(c) [below right = 4mm  of 3nn]{(c)};

\end{scope}
    \begin{scope}
      \draw (1) -- (2);
      \draw (1) -- (3);
	\draw (2) -- (4);
	\draw (3) -- (4);
    \end{scope}
    \begin{scope}[->, > = latex']
    \draw (1n) -- (2n);
    \draw (1n) -- (3n);
	\draw (2n) -- (4n);
	\draw (3n) -- (4n);
    \draw (1nn) -- (2nn);
    \draw (1nn) -- (3nn);
	\draw (2nn) -- (4nn);
    \end{scope}
    \begin{scope}[<->, > = latex']
    \draw (3nn) -- (4nn);
    \end{scope}

    \end{tikzpicture}
		\caption{{\footnotesize (a) The skeleton of $P$ with no arrows. (b) A DAG $G_1$.  (c) A graph  $G_2$ with  arcs.}}\label{fig:exunique}
		\end{figure}
\end{eexample}

From  Example \ref{ex:10}, we see that  singleton-transitivity is not necessary for the uniqueness property.  We provide a (weaker) condition, which together with assumptions of Theorem \ref{thm:ske},  are still sufficient for uniqueness.

We say that an independence model $\mathcal{J}$ is \emph{path-unstable} if there exists a  path
$$\langle i=i_0,i_1,\cdots,i_r,k,j\rangle, \nn\nn i\not\sim j \text{ and }  i_r\sim j \text{ for every } r \geq 1$$ in $\sk(\mathcal{J})$ such that $i\ci j\cd U\cup C$ and $i\ci j\cd U\cup C\cup\{k\}$, for some $C$  disjoint from $\{i,j,k\}$, where  $U=\{i_1,\cdots,i_r\}$ or $U = \emptyset$ when $r =0$. An independence model that is not path-unstable is called \emph{path-stable}.

In the case where the structure-learning algorithm is DAG-generating, the required condition is simpler. We say that an independence model $\mathcal{J}$ is \emph{V-unstable} if there exists a  V-configuration  $\langle i,k,j\rangle$ in $\sk(\mathcal{J})$ such that $i\ci j\cd C$ and $i\ci j\cd C\cup\{k\}$, for some $C$  disjoint from $\{i,j,k\}$. An independence model that is not V-unstable is called \emph{V-stable}. %
In the same fashion as before, if $\mathcal{J}(P)$ satisfies path- or V-stability then we simply say that ``$P$ satisfies these conditions."

\begin{remark}
As will be seen in the proof of Proposition \ref{thm:unique-sin-up-st}, the path that appears in the definition of a path-unstable $P$ can be thought of as the skeleton of (and can be directed to become) a \emph{discriminating path} \citep{ali09} from $i$ to $j$; i.e., a path $\langle i=i_0,i_1,\cdots,i_r,k,j\rangle$ that is a collider path, and  $i_n\in\pa(j)$, for every $1\leq n\leq r$.
\end{remark}
\begin{remark}
The well-known
{orientation-faithfulness}
assumption (which is defined as a part of the
{restricted-faithfulness}
assumptions) \citep{ram06} states that for all V-configurations $\langle j, l, k\rangle$ and all subsets $S \subset V \setminus \{j, k\}$ such that $j$ is d-connected to $k$ given $S$ it holds that $j \notci k \cd S$. Notice that, unlike orientation-faithfulness,  V-stability (and also path stability) are purely defined on the distribution $P$ and not on the unknown causal graph $G_0$.  Indeed, if the conditions of Theorem \ref{thm:ske} are satisfied to ensure that $\sk(G_0)=\sk(G(P))$,
then orientation-faithfulness implies V-stability.  However, even under these extra assumptions, the other direction does not hold; see Example \ref{ex:or-faith} below.
\end{remark}
\begin{eexample}\label{ex:or-faith}
Consider the true causal graph $G_0$ in Figure \ref{fig:or-faith}. Let $\mathcal{J}$ contain all the conditional independence statements implied by the global Markov property on $G_0$,  with the addition of  $j\ci k\cd s$.  We have that
$\sk(\mathcal{J})=\sk(G_0)$.  By looking at
$\sk(\mathcal{J})$,
it can be seen that V-stability is satisfied, but orientation-faithfulness is not: Consider the V-configuration $\langle j, l,k\rangle$.  We have that $j$ is d-connected to $k$ given $s$,  but that $j \ci k \cd s$.
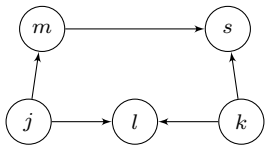
\begin{figure}[htb]
\centering
\begin{tikzpicture}[node distance = 8mm and 8mm, minimum width = 6mm]
    \begin{scope}
      \tikzstyle{every node} = [shape = circle,
      font = \scriptsize,
      minimum height = 6mm,
      inner sep = 2pt,
      draw = black,
      fill = white,
      anchor = center],
      text centered]
      \node(1) at (0,0) {$l$};
      \node(2) [above left = of 1] {$m$};
      \node(3) [above right = of 1] {$s$};
      \node(4) [right = of 1] {$k$};
      \node(5) [left = of 1] {$j$};
    \end{scope}
		
    \begin{scope}[->, > = latex']
    \draw (5) -- (1);
    \draw (4) -- (1);
	\draw (5) -- (2);
	\draw (4) -- (3);
    \draw (2) -- (3);
    \end{scope}

\end{tikzpicture}
		\caption{{\footnotesize A directed acyclic graph $G_0$.}}\label{fig:or-faith}
		\end{figure}

\end{eexample}
V-stability is  weaker than singleton-transitivity:
\begin{proposition}
   Singleton-transitivity implies V-stability.
\end{proposition}
\begin{proof}
Towards a contradiction, suppose that there is a V-configuration $\langle i,k,j\rangle$ in $\sk(P)$ such that $i\ci j\cd C$ and $i\ci j\cd C\cup\{k\}$, for some $C$  disjoint from $\{i,j,k\}$. By singleton-transitivity, $i\ci k\cd C$ or $k\ci j\cd C$,  which is absurd since both $ik$ and $jk$ are edges in $\sk(P)$.
\end{proof}
In addition, the converse of the above result does not hold; see Example \ref{ex:10}. Hence, we have the following corollary:
\begin{corollary}
\label{logically-implied}
If $P$ is faithful to a graph, then it satisfies V-stability.
\end{corollary}
\begin{proof}
The proof follows from the fact that $P$ being graphical implies it satisfies singleton-transitivity (Proposition \ref{thm:2}).
\end{proof}
Also notice that singleton-transitivity and path-stability do not imply each other (for one direction, see Example \ref{ex:10} again), but there are considerably fewer independence statements to test for for path-stability than for singleton-transitivity. In a similar fashion, there could be cases where $P$ is graphical, but path-stability is not satisfied: one can take the skeleton of a discriminating path and direct the edges (with the path not being collider) so that the induced independence model (by separation) is path unstable.
However, notice again that unlike faithfulness, path-stability is purely defined on $P$, and hence, in principle, testable. We now have the following relationship to the uniqueness property:
\begin{proposition}
	\label{thm:unique-sin-up-st}
	Let $P$ be a distribution.
\begin{enumerate}[(a)]
  \item
  \label{path-U}
  If $P$ is path-stable, then it satisfies the uniqueness property.
  \item
  \label{DAG-U}
  If $P$ is V-stable, then it satisfies the DAG-uniqueness property.
\end{enumerate}
\end{proposition}
\begin{proof}
We start with the proof for the case of DAGs, which is more elementary than the general case.

(b) Suppose that $P$ satisfies ordered upward- and downward-stability w.r.t.\ two orders $\leq'$ and $\leq''$.  Towards a  contradiction, suppose that the two DAGs,  $G(P,\leq')$ and $G(P,\leq'')$ are not Markov equivalent.
	
	 By Lemma \ref{lem:Mark-eq-DAG},  there is a collider V-configuration $\langle i,k,j\rangle$, say, in $G(P,\leq')$, which is a non-collider in $G(P,\leq'')$ (say $k>''j$). 
	Since  $G(P,\leq')$ is maximal, $i\ci j\cd C$, for some $C$.   By ordered downward-stability w.r.t.\ $\leq'$, we obtain $i\ci j\cd C'$, for some $C'$ such that $k\notin C'$. Now by ordered upward-stability w.r.t.\ $\leq''$, we obtain $i\ci j\cd C'\cup\{k\}$. This shows that $P$ is V-unstable, which is a contradiction.

(a) Suppose that $P$ satisfies ordered upward- and downward-stability w.r.t.\ two orders $\leq'$ and $\leq''$.  Towards a  contradiction, suppose that the two maximal graphs $G(P,\leq')$ and $G(P,\leq'')$ are not Markov equivalent.

By Lemma \ref{lem:Mark-eq-MAG}, there is a minimal collider path $\pi$, say, in $G(P,\leq')$, which is not a minimal collider path in $G(P,\leq'')$. 
Consider a shortest minimal collider subpath $\pi'$ of $\pi$ in $G(P,\leq')$ that is  not a  minimal collider in $G(P,\leq'')$.   Furthermore, we can assume that $\pi'$ is not a collider path in $G(P,\leq'')$, since if
 $\pi'$ is a collider path, but not minimal, we  consider the minimal collider subpath $\pi''$ in $G(P,\leq'')$.  Note that $\pi''$ is not a minimal collider in $G(P,\leq')$. We can then continue the proof using $\pi''$ instead of $\pi'$, and the roles of $G(P,\leq')$ and $G(P,\leq'')$ reversed.

Let $i$ and $j$ be the endpoints of $\pi'$ which is not a collider in $G(P,\leq'')$.  Consider a node $k$, which is a non-collider node on $\pi'$.  We show that $k$ is adjacent to $i$ or $j$ on $\pi'$.
Towards a contradiction, suppose that it is not adjacent to either.
Consider a  subpath of $\pi'$ given by $\langle i',h,k,\ell,j'\rangle$.  Recall that the endpoints of a minimal collider path are not adjacent.  Hence $h\sim \ell$, otherwise, $\langle h,k,\ell\rangle$ is a shorter minimal collider subpath in $G(P,\leq')$ that is not a minimal collider in $G(P,\leq'')$.   Similarly, $h\sim j'$ and $i'\sim \ell$. Since $\pi'$ is a minimal collider, the $h\ell$-edge cannot be an arc, and the $hj$- and $i\ell$-edge have to be arrows from $h$ to $j'$ and $\ell$ to $i'$, respectively.  Since directed cycles are not permissible in an ancestral graph, it is easy to rule out the existence of an arrow from $h$ to $\ell$ and vice-versa; we already ruled out the possibility that the $h\ell$-edge is an arc, and we require that $h \sim \ell$, which is absurd.

 Without loss of generality, assume that $k$ is adjacent to $j$ on $\pi'$.
Let $\pi'=\langle i=i_0,i_1,\cdots,i_r,k,j\rangle$, where $r$ may be $0$. If $r>0$, then by minimality, $i_r\sim j$ and $i_r \fra j$.   From an  inductive argument along $\pi'$, if $r>0$,  we conclude that $\pi'$ is a discriminating path from $i$ to $j$.

Since $G(P,\leq')$ is maximal, $i\ci j\cd C$, for some $C$; moreover,  $k\notin\an(i)$ since otherwise $i$ and $j$ are not separated.  Since $k$ is a collider in $G(P, <')$, we have $k \not \in an(j)$.   Therefore, by ordered downward-stability w.r.t.\ $\leq'$, we obtain $i\ci j\cd C'$, for some $C'$ such that $k\notin C'$. By ordered upward-stability w.r.t.\ $\leq'$, and the previously established discrimination,  we obtain $i\ci j\cd U\cup C'$, where $U=\{i_1,\cdots,i_r\}$ (and if $r =0$, then $U = \emptyset$).    Now by ordered upward-stability w.r.t.\ $\leq''$, we obtain $i\ci j\cd U\cup C'\cup\{k\}$, since $k$ is not a collider node on $G(P, \leq'')$.   Thus $P$ is path-unstable, which is a contradiction.
\end{proof}

This leads to the main result, which provides conditions under which
a natural
 structure-learning algorithm works.
\begin{theorem}\label{coro:main}
Suppose that $P$ satisfies the conditions of Theorem \ref{thm:ske}, that is, it is Markovian and converse pairwise Markovian to $G_0$, and satisfies ordered upward- and downward-stability w.r.t.\ $G_0$.
Let $G(P)$ be the output of a natural structure-learning algorithm.
\begin{enumerate}[(a)]
  \item
  If  $P$ is also path-stable, then  $G(P)$ is Markov equivalent to $G_0$.
  \item
  If $G_0$ is a DAG and $P$ is V-stable, then a DAG-output $G(P)$ is Markov equivalent to $G_0$.
\end{enumerate}
\end{theorem}
\begin{proof}
The proof follows from Proposition \ref{thm:unique-sin-up-st} and Lemma \ref{lem:essen}.
\end{proof}
As an example for the theorem above,  we can revisit our motivating Example \ref{ex:10} to observe that although faithfulness is not satisfied, all the conditions of Theorem \ref{coro:main}, and, in particular, V-stability are satisfied.

However, notice that, under the conditions presented here, ordered downward- and upward-stability w.r.t.\ $G_0$ are not necessary in the sense that there are cases where ordered downward- and upward-stability  are not satisfied w.r.t.\ $G_0$ , but an output of the algorithm is Markov equivalent to $G_0$. This stems from the fact that there could be two Markov equivalent graphs and a $P$ such that $P$ satisfies ordered downward- and upward-stability w.r.t.\ one but not the other; see example below.
\begin{eexample}
Consider the causal graph $G_0$ and its Markov equivalent graph $G_1$ below. Let $\mathcal{J}$ (which, in principle could be induced by a probability distribution) contain all the independencies induced by the global Markov property plus $k\ci m\cd l$ and additional independence statements implied by applying ordered upward-stability w.r.t.\ $G_1$ and semi-graphoid axioms. It can be checked that $k\notci m\cd \{i,l\}$. Hence, $\mathcal{J}$ does not satisfy ordered upward-stability w.r.t.\ $G_0$. However, by definition, $\mathcal{J}$ satisfies ordered upward-stability w.r.t.\ $G_1$. Notice also that, in this example, $V$-stability (and path-stability) are satisfied.
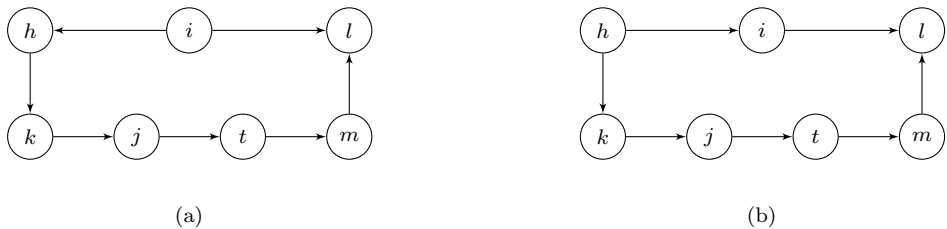
\begin{figure}[htb]
\centering
\begin{tikzpicture}[node distance = 8mm and 8mm, minimum width = 6mm]
    \begin{scope}
      \tikzstyle{every node} = [shape = circle,
      font = \scriptsize,
      minimum height = 6mm,
      inner sep = 2pt,
      draw = black,
      fill = white,
      anchor = center],
      text centered]
      \node(1) at (0,0) {$j$};
      \node(2) [left = of 1] {$k$};
      \node(3) [above  = of 2] {$h$};
      \node(7) [right  = of 1] {$t$};
      \node(4) [right = of 7] {$m$};
      \node(5) [above  = of 4] {$l$};
      \node(6) [right = 15mm of 3] {$i$};
      \node(1n) [right = 70mm of 1] {$j$};
      \node(2n) [left = of 1n] {$k$};
      \node(3n) [above  = of 2n] {$h$};
      \node(7n) [right  = of 1n] {$t$};
      \node(4n) [right = of 7n] {$m$};
      \node(5n) [above = of 4n] {$l$};
      \node(6n) [right = 15mm of 3n] {$i$};
    \end{scope}
		
    \begin{scope}
    \tikzstyle{every node} = [node distance = 6mm and 6mm, minimum width = 6mm,
    font= \scriptsize,
      anchor = center,
      text centered]
\node(a) [below = 19mm of 6]{(a)};
\node(b) [below = 19mm  of 6n]{(b)};

\end{scope}

    \begin{scope}[->, > = latex']
    \draw (6) -- (5);
    \draw (6) -- (3);
	\draw (4) -- (5);
	\draw (3) -- (2);
    \draw (2) -- (1);
    \draw (1) -- (7);
    \draw (7) -- (4);
    \draw (6n) -- (5n);
    \draw (3n) -- (6n);
	\draw (4n) -- (5n);
	\draw (3n) -- (2n);
    \draw (2n) -- (1n);
    \draw (1n) -- (7n);
    \draw (7n) -- (4n);
    \end{scope}

    \end{tikzpicture}
		\caption{{\footnotesize (a) $G_0$, (b) $G_1$}}\label{fig:ex0}
		\end{figure}
\end{eexample}
The above discussion suggests that one can replace the condition ``$P$ satisfying ordered upward- and downward-stability w.r.t.\ $G_0$" in Theorem \ref{coro:main} with ``$P$ satisfying ordered upward- and downward-stability w.r.t.\ a graph in the Markov equivalence class of $G_0$" in order to make the condition necessary.
\section{Structure learning for structural causal models}\label{sec:st-l-SCM}

First, we give a proof for  the global Markov property for general SCMs, which allow for the possibility of arcs.  Then we find conditions on the functions of the SCM to ensure the converse pairwise Markov property; see \eqref{converse-def}.

\subsection{Global Markov property for general SCMs}
The main assumption used in the general setting for structure learning is that $P$ is Markovian to the causal graph $G_0$. A main reason for popularity of SCMs for causal learning is that this condition is automatically satisfied for the joint distribution of an SCM. Indeed, the joint distribution of an SCM with mutually independent errors being Markovian to the DAG $G_0$ is a well-known result \citep{ver88,pea09}, but we
are
not aware of such a result in an explicit and non-parametric form when unobserved variables are existent; see \citet{kos99} for a similar result in the Gaussian case and with a slightly different type of graph. Recently, in \citet{bon20}, a similar result is proven for directed SCMs, but for a less general setting when one focuses on ancestral graphs; see Remark \ref{rem:joint-noise}.  See also \cite{forre2017markov} for various equivalent notions of the Markov property.
\begin{theorem}
\label{prop:scm-mark}
The joint distribution $P$ of an SCM is Markovian to\ $G_0$.
\end{theorem}

The proof of the above theorem requires some definitions and results that are not used in the  subsequent sections.  We provide the proof together with these definitions and results in the remainder of this subsection. Readers may skip the proof without loss of continuity.

We first need some primary definitions and results mostly from \citet{ric03}. The \emph{district} of node $i$, denoted by $\dis(i)$, is the set of nodes connected to $i$ by a collider path on which every edge is an arc. We also define a set $A$ to be \emph{ancestral} if it is closed under the ancestor relation, i.e., if $A\cup\an(A) = A$. If $A$ is an ancestral set in a graph $G$, and $i$ is a node in $A$ that has no children in $A$, then we define the \emph{Markov blanket} of $i$ w.r.t.\ the induced subgraph $G[A]$ on $A$ by $$\mb(i,A)=\pa_{G[A]}(\dis_{G[A]}(i)) \ \cup \    (\dis_{G[A]}(i)\setminus\{i\}).$$

We use the notation $\pst(i)$ for the set of nodes with an order larger that $i$. We say that a probability distribution satisfies the \emph{ordered local Markov property} w.r.t.\ the ordering $\leq$, if for any node $i$, and ancestral set $A$ such that $i\in A\subseteq \pst(i)$, we have
$$i\ci  [ A\setminus (\mb(i,A)\cup \{i\}) ] \cd \mb(i,A).$$
We have the following equivalence, which does not require any extra assumptions, originally proven for ADMGs (which is a more general class than directed ancestral graphs):
\begin{lemma}[\cite{ric03}]
	\label{lem:ord-mar-eq}
For a valid ordering $\leq$ of the nodes of a graph $G$, a distribution $P$ satisfying the ordered local Markov property w.r.t.\ $\leq$ is equivalent to $P$ being Markovian to $G$.
\end{lemma}
In this section, we also utilize the concept of marginalization over probability distributions and graphs. Let $\alpha(P,M)$ be the marginalized distribution over $M\subset V$ so that $\alpha(P,M)$ is defined over $X_{V\setminus M}$. Notice that the induced independence model of $\alpha(P,M)$ is $\mathcal{J}(\alpha(P,M))=\{A\ci B\cd D:(A\cup B\cup D)\cap M=\varnothing\}$.

By marginalizing over the node set $M$ of graph $G$, a new graph $\alpha(G,M)$ is obtained from which the node set $M$ is removed: For every path $\pi$ in $G$ between $i,j\notin M$, whose inner nodes are all non-collider and in $M$, we generate an edge in $\alpha(G,M)$, which preserves the (lack of) arrowheads pointing to $i$ and $j$ in $\pi$; We apply this repeatedly until no new edges are being generated, and then remove all the nodes in $M$; for more details of this algorithm, refer to \citet{ric02,sad13}. The following result is well-known \citep{ric02,sad13}:
\begin{lemma}\label{lem:marg0}
For a graph $G$ and disjoint node sets $A$, $B$, $C$, and $M$, it holds that
\begin{displaymath}
A\dse B\cd C\text{ in }\alpha(G,M) \iff
A\dse B\cd C\text{ in }G.
\end{displaymath}
\end{lemma}
We shall need the following corollary of this result:
\begin{lemma}\label{lem:marg}
Suppose that a distribution $P$ is Markovian to the graph $G$. Then $\alpha(P,M)$ is Markovian to $\alpha(G,M)$.
\end{lemma}
\begin{proof}
$P$ being Markovian to $G$ implies that $A\dse B\cd C$ in $G$ implies that $A\ci B\cd C$. $A\dse B\cd C$ in $\alpha(G,M)$ implies that $(A\cup B\cup C)\cap M=\varnothing$, and also $A\dse B\cd C$ (by Lemma \ref{lem:marg0}).  $(A\cup B\cup C)\cap M=\varnothing$ implies that  $A\ci B\cd C$ exists in $\mathcal{J}(\alpha(P,M))$. Therefore, $A\dse B\cd C$ in $\alpha(G,M)$ implies $A\ci B\cd C$ exists in $\mathcal{J}(\alpha(P,M))$.
\end{proof}

	We are now ready to provide the proof of Theorem \ref{prop:scm-mark}.

\begin{proof}{\bf of Theorem \ref{prop:scm-mark} \ }
Consider the augmented graph $G_0(X,\epsilon)$ and the joint distribution $P_{X,\epsilon}$. Consider the ordering similar to the ordering in the SCM for $X_i$s, and, for $\epsilon_i$s, let $X_{i+1}>\epsilon_i>X_i$. We show that $P_{X,\epsilon}$ satisfies the ordered local Markov property w.r.t.\ this ordering in $G_0(X,\epsilon)$. Hence by Lemma \ref{lem:ord-mar-eq}, we obtain that $P_{X,\epsilon}$ is Markovian to  $G_0(X,\epsilon)$.  Consider a node in $G_0(X,\epsilon)$. We have two cases:

First, assume that this node is of $\epsilon$ form, i.e., say it is $\epsilon_i$. Consider also $A\subseteq \pst(\epsilon_i)$. We have that $\mb(\epsilon_i,A)=\{\epsilon_j\in A:\n X_j\in\dis(X_i) \text{ in } G_0\}$. If $\mb(\epsilon_i,A)=A$, then $\epsilon_i\ci [A\setminus (\mb(\epsilon_i,A)\cup\{\epsilon_i\})]   \cd \mb(\epsilon_i,A)$ is trivial; otherwise there is no edge between $\mb(\epsilon_i,A)$ and $A\setminus (\mb(\epsilon_i,A)\cup\{\epsilon_i\})$ in  $G_0(X,\epsilon)$, which, by joint independence of such $\epsilon_j$, implies that $\{\epsilon_i\}\cup\mb(\epsilon_i,A)\ci A\setminus (\mb(\epsilon_i,A)\cup\{\epsilon_i\})$. This, by the weak union property, implies $$\epsilon_i\ci [A\setminus (\mb(\epsilon_i,A)\cup\{\epsilon_i\})]\cd \mb(\epsilon_i,A).$$

Secondly, assume that this node is of $X$ form, i.e., say it is $X_i$.  Consider also $A\subseteq \pst(X_i)=(\epsilon_i,X_{i+1},\epsilon_{i+1},\dots,X_n,\epsilon_n)$. We have that $\mb(X_i,A)=\{\epsilon_i\}\cup\{X_j:\n X_j\in\pa_{G_0[A]}(X_i)\}$. However, given $\mb(X_i,A)$, we have that $X_i$ is deterministic; hence, it is independent of all disjoint random vectors including $A\setminus (\mb(X_i,A)\cup\{X_i\})$.

Finally, since the inner nodes of $\langle X_i,\epsilon_i,\epsilon_j,X_j\rangle$ are non-collider in $\alpha(G_0(X,\epsilon)$, it can be seen that $\alpha(G_0(X,\epsilon),\epsilon)=G_0$. Hence, $P_{X,\epsilon}$ being Markovian to  $G_0(X,\epsilon)$, by Lemma \ref{lem:marg}, implies that $P$ is Markovian to $G_0$.
\end{proof}
\subsection{Conditions and assumptions for structure learning of SCMs}
The general setting for structure learning for SCMs is the same as in the previous sections, and, in fact, the case of SCM could be considered a special case of the theory presented above: We start by assuming the existence of an \emph{unknown} SCM $\mathfrak{C}$ associated to an \emph{unknown} causal graph $G_0$ as described in Section \ref{sec:SCM}. For structure learning (from observational data), we find a graph that belongs to the Markov equivalence class of $G_0$. As discussed, the SCM $\mathfrak{C}$ entails a joint distribution $P$, which is still unknown. Therefore, the relationships can be illustrated in the following diagram:
$$G_0 \leftarrow\mathfrak{C}\rightarrow P \stackrel{\text{conditional independence}}{\rightarrow} \mathcal{J}(P)\stackrel{\text{Algorithm}}{\rightarrow} G(P) \stackrel{\text{Markov equivalent}}{\sim} G_0$$

We begin our discussion with a few examples.
In general, the independence models of the SCMs are not faithful to $G_0$.  Although there are many examples where faithfulness fails, the following simple and perhaps familiar  example will help to motivate further considerations regarding the converse pairwise Markov property.
\begin{eexample}[Failure of faithfulness]
	\label{mod}
	Consider the following SCM where $2 \fra 1$, and $2$ is the only parent of $1$.    Suppose $X_2 = \epsilon_2$ and $X_1 = X_2 \oplus \epsilon_1 \bmod 2$, where $\epsilon_2$ and $\epsilon_1$ are independent Bernoulli random variables with parameter $p=\tfrac{1}{2}$.  It is easy to verify that $X_2$ is also independent of $X_1$.
	
	For a continuous example, we exploit of the invariance of Lebesgue measure on the circle.   Take $\epsilon_1$ and $X_2=\epsilon_2$ to be independent continuous real-valued random variables.  Let $F$ and $G$ be the cumulative distribution functions for $X_2$ and $\epsilon_1$, respectively.  Set  $X_1 = F(X_2) \oplus G(\epsilon_1) \bmod 1$.   It is easy to verify that $X_2$  is independent of $X_1$.
		\end{eexample}

 Motivated by Example \ref{mod} and the elementary fact that  a random rotation that is independent of a standard bivariate normal remains independent of the randomly rotated bivariate normal, and more generally the connection between symmetries and independence \citep{MR813238}, we have the following lemma.

\begin{lemma}
\label{char}
Consider the SCM with associated graph $2 \fra 1$, so that  $2$ is the only parent of $1$.  Suppose that $X_2 = \epsilon_2$, $X_1 = \phi_1(X_2, \epsilon_1)$ and that $\epsilon_2$ and $\epsilon_1$ are independent.   Then
$X_2$ is independent of $X_1$ if and only
\begin{equation}
	\label{dist-e}
\phi_1(x_2, \epsilon_1) \stackrel{d}{=} X_1
\end{equation}
for almost every $x_2 \in \mathcal{X}_2$.
	\end{lemma}

Note that \eqref{dist-e} is easily satisfied if $\phi_1$ is actually a function of $\epsilon_1$ alone. In \eqref{dist-e}, it may be helpful to think of $x_2$ as indexing (distinct) ways of generating a random variable with law $X_1$, given the randomization $\epsilon_1$.

\begin{proof}{\bf of Lemma \ref{char} \ }
Let $S_2 \subset  \mathcal{X}_2$ and $S_1 \subset  \mathcal{X}_2$ be measurable subsets.  Let $Q_{X_2}$ and $Q_{\epsilon_1}$ be the laws of $X_2$ and $\epsilon_1$, respectively.   Since $2 \fra 1$, and  $X_2$ and $\epsilon_1$ are independent, we have
		\begin{eqnarray*}
			\mathbb{P}(X_1 \in S_1, X_2 \in S_2) &=& \mathbb{P}(\phi_1(X_2, \epsilon_1) \in S_1, X_2 \in S_2) \\
			&=& \int \int \mathbf{1} [\phi_1(x_2,e_1) \in S_1)] \mathbf{1}[x_2 \in S_2] dQ_{X_2}(x_2) d Q_{\epsilon_1}(e_1)  \\
			&=& \int_{S_2}  \mathbb{P} (\phi_1(x_2,\epsilon_1) \in S_1) dQ_{X_2}(x_2).
		\end{eqnarray*}
Thus if \eqref{dist-e} holds, we have
\begin{eqnarray*}
	\mathbb{P}(X_1 \in S_1, X_2 \in S_2) &=& \int_{S_2} \mathbb{P}(X_1 \in S_1) dQ_{X_2}(x_2) \\
	&=&
	 \mathbb{P}(X_1 \in S_1) Q_{X_2}( S_2)\\
	 &=&\mathbb{P}(X_1 \in S_1) \mathbb{P}( X_2 \in S_2)
	\end{eqnarray*}
	and the desired independence follows.

If $X_2$ and $X_1$ are independent, then
$$	\mathbb{P}(X_1 \in S_1) =
	 \frac{1}{ \mathbb{P}(X_2 \in S_2)}\int_{S_2}  \mathbb{P}( \phi_1(x_2,\epsilon_1)) \in S_1) dQ_{X_2}(x_2)$$
	 for all measurable subsets $S_2$ with $Q_{X_2}(S_2) >0$.   By setting $$S_2 =  \{ x_2 \in \mathcal{X}_2 : \mathbb{P}(\phi_1(x_2, \epsilon_1) \in S_1) > \mathbb{P} (X_1 \in S) \},$$
	 we see that $Q_{X_2}(S_2) = 0$.  Hence together with  reversing the inequality in the definition of $S_2$, we obtain  \eqref{dist-e}.
	\end{proof}

In fact, the joint distributions of SCMs may not satisfy any of the axioms for faithfulness, presented in Proposition \ref{thm:2}, and discussed in the Sections \ref{sec:propind} and \ref{sec:ouds}. For a counter-example  and a thorough discussion on the intersection property, see \citet{pet15}. For an example of an SCM not satisfying singleton-transitivity, see Example \ref{ex:1}. For counter-examples for other axioms, we provide the following examples:
\begin{eexample}[Failure of composition and ordered upward-stability]
Consider the following SCM associated to the graph $G_0$ in Figure \ref{fig:excomp}:
$$X_3=\epsilon_3;\nn X_2=\epsilon_2;\nn X_1=X_2\oplus X_3\oplus \epsilon_1 \bmod 2;$$
where $\epsilon_1,\epsilon_2,$ and $\epsilon_3$ are independent Bernoulli random variables, $\epsilon_2$ and  $\epsilon_3$ with parameter $\tfrac{1}{2}$, and $\epsilon_1$ with parameter $\tfrac{1}{3}$.

It is easy to verify that $X_1\ci X_2$, and $X_1\ci X_3$, but $X_1\notci (X_2,X_3)$ and $X_1 \notci X_2\cd X_3$. Hence this  is an example for an SCM not satisfying the composition property or ordered upward-stability.
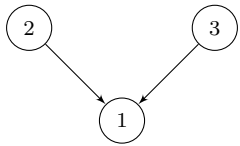
\begin{figure}[htb]
\centering
\begin{tikzpicture}[node distance = 8mm and 8mm, minimum width = 6mm]
    \begin{scope}
      \tikzstyle{every node} = [shape = circle,
      font = \scriptsize,
      minimum height = 6mm,
      inner sep = 2pt,
      draw = black,
      fill = white,
      anchor = center],
      text centered]
      \node(1) at (0,0) {$1$};
      \node(2) [above left = of 1] {$2$};
      \node(3) [above right = of 1] {$3$};

    \end{scope}

    \begin{scope}[->, > = latex']
    \draw (2) -- (1);
    \draw (3) -- (1);

    \end{scope}

\end{tikzpicture}
		\caption{{\footnotesize A directed acyclic graph $G_0$.}}\label{fig:excomp}
		\end{figure}

\end{eexample}

\begin{eexample}[Failure of ordered downward-stability]
Consider a complete DAG with nodes $(1,2,3)$. Let $\epsilon_3 = (N_3, U_3)$ and $\epsilon_2 = (U_2, V_2)$, where $N_3, U_3, U_2$, and $V_2$ are all independent, and $N_3$ is Poisson distributed, and $U_3, U_2$ and $V_2$ are uniformly distributed on $[0,1]$.   Let $X_3 = \epsilon_3$.    Let $X_2=\phi_2(X_3, \epsilon_2) = (\psi_2(N_3, U_2), V_2)$, where $\psi_2$ is a deterministic function such that $\psi_2(n_3, U_2)$ is uniformly distributed on the integers $\{1, \ldots, n_3+1\}$.     Let $X_1=\phi_1(X_3, X_2) = (N_3, V_2)$.  Clearly, $X_2$ is not independent of $X_3$.  However, conditional on $X_1$, we have  that $X_2$ is independent of $X_3$.
\end{eexample}

An important property used in the results of previous section is converse pairwise Markov property. Here, we are mainly concerned with providing conditions for verifying this property. For an SCM $\mathcal{C}$, we provide an assumption on the joint distribution of $\mathcal{C}$ that implies the converse pairwise Markov property.  Consider the following extension of Lemma \ref{char}.

\begin{proposition}
	\label{pairwise-distinct}
	  Consider an  SCM, where $j \fra i$.  Write  $X_i= \phi_i(X_{\pa(i) \setminus \{j\} }, X_j, \epsilon_i)$.  Suppose that $\an(i,j) \not =\emptyset$.    Then   $i\ci j\cd \an(i,j)$ if and only if for almost all $x_{\an(i,j)} \in \mathcal{X}_{\an(i,j)}$ the conditional law of $X_i$ given $X_{\an(i,j)} = x_{\an(i,j)}$ is the same as the law of $\phi_i(x_{\pa(i) \setminus \{j\}}, x_j, \epsilon_i)$
	for almost all $x_j \in \mathcal{X}_j$ in the support of the conditional distribution of $X_j$ given $x_{\an(i,j)}$.
	\end{proposition}

\begin{proof}  Let $\mathcal{F}_j$ and $\mathcal{F}_i$ be the $\sigma$-algebras for $\mathcal{X}_j$ and $\mathcal{X}_i$.   Let $D:\mathcal{X}_{\an(i,j)} \times (\mathcal{F}_j, \mathcal{F}_i) \to [0,1]$ be a  regular conditional distribution for $(X_j, X_i)$ given $X_{\an(i,j)}$.  Fix $x_{\an(i,j)} \in \mathcal{X}_{\an(i,j)}$.  All the parents of $i$ are included in $\an(i,j)$, except $j$.    Let $Q_i$ and $Q_j$ be marginal distributions of $X_i$ and $X_j$ under $D_{x_\an(i,j)} = D(x_\an(i,j), \cdot, \cdot)$.
	Following the proof of Lemma \ref{char}, since $j \fra i$ and the graph is ancestral,  we have that  $\epsilon_i$, with law $Q_{\epsilon_i}$,   is independent of $(X_{\an(i,j)}, X_j)$.  Thus  for $S_j \in \mathcal{F}_j$ and $S_i \in \mathcal{F}_i$, we have
		\begin{eqnarray*}
			D_{x_{\an(i,j)}}(S_j, S_i)
			&=& \int \int_{S_j} \mathbf{1} [\phi_i(x_{{\pa(i) \setminus \{j\} }},x_j,e_i) \in S_i)]  dQ_j(x_j) d Q_{\epsilon_i}(e_i)  \\
			&=& \int_{S_j}  \mathbb{P} (\phi_i(x_{{\pa(i) \setminus \{j\} })},x_j,\epsilon_i) \in S_i)  dQ_j(x_j).
		\end{eqnarray*}
Thus, if
\begin{equation}
	\label{dist-eq}
 Q_i( S_i) = 	\mathbb{P} (\phi_i(x_{{\pa(i) \setminus \{j\} })},x_j,\epsilon_i) \in S_i)
 \end{equation}
for all $Q_j$-almost all $x_j \in \mathcal{X}_j$, then
$$	D_{x_{\an(i,j)}}(S_j, S_i) = Q_j(S_j) Q_i( S_i)$$
as desired.

A routine variation of the rest of the proof of Lemma \ref{char}, gives that the condition on $\phi_i$ expressed in \eqref{dist-eq} is also necessary for the conditional independence of $X_j$ and $X_i$.
	\end{proof}

Proposition \ref{pairwise-distinct} includes Lemma \ref{char} by making the suitable interpretations in the case that $\an(i,j) = \emptyset$.

   Proposition \ref{pairwise-distinct} suggests the following condition that implies the converse pairwise Markov property, in the case that true causal graph is a DAG.
We say that the (joint distribution) of an  SCM satisfies \emph{positivity} if   $\mathbb{P}(X_B \in S|X_A = x_A) >0$ almost surely, for every set of nodes $B$ disjoint from $A$ and every measurable $S$ for which $\mathbb{P}(X_B \in S) >0$, for almost every $x_A \in \mathcal{X}_A$.

Observe that each function $\phi$ in an SCM has the form  $$\phi:\mathcal{X}_1 \times \cdots \times \mathcal{X}_n \times \mathcal{E} \to \mathcal{X}_{n+1},$$ which is endowed with a  product probability measure $\mu \otimes Q_{\epsilon}$, the joint law of the parents and exogenous variable on  $\mathcal{X}_1 \times \cdots \times \mathcal{X}_n \times \mathcal{E}$, and a probability measure $\nu = (\mu \otimes Q_{\epsilon}) \circ \phi^{-1}$ on   $\mathcal{X}_{n+1}$, the law of the child.  Let $\mu^j$ be the projection of $\mu$ to space
$$\mathcal{X}^j := \mathcal{X}_1 \times \cdots \times \mathcal{X}_{j-1} \times \mathcal{X}_{j+1} \times \cdots \times \mathcal{X}_n.$$
   Let $\epsilon$, $X= (X_1, \ldots, X_n)$, and $X_j$ have laws $Q_{\epsilon}$, $\mu$, and $\mu_j$ respectively. We say that $\phi$ has  \emph{non-constant fibers} if for
   all
    $1 \leq j \leq n$
    there is some
     $F \subset \mathcal{X}^j$ with positive $\mu^j$-measure  for each
     $$x^j=(x_1, \ldots, x_{j-1}, x_{j+1}, \ldots, x_n) \in F$$ and there is some  $S \subset \mathcal{X}_{n+1}$ such that  the function $k: \mathcal{X}_j \to [0,1]$ given by

$$k(x_j)= \mathbb{P}(\phi(x_1, \ldots, x_j, \ldots, x_n, \epsilon) \in S)$$
is non-constant $\mu_j$-almost surely.   We say that an SCM has   \emph{non-constant fibers} if every function in the SCM has non-constant fibers.
\begin{remark}
	\label{fibers}
Note that in the proof of Proposition \ref{pairwise-distinct}, equality \eqref{dist-eq} fails if  $\phi_i$ has non-constant fibers and the SCM satisfies positivity.
\end{remark}

\begin{remark}
	\label{additive}
	It is easy to see that additive noise models with the form $$\phi(x, \epsilon) = h(x) + \epsilon$$ will have non-constant fibers
	under the mild condition that for any $1 \leq j \leq n$, and any fixed $x_1, \ldots, x_{j-1}, x_{j+1}, \ldots, x_n$, the function
	$$g(x_j) = h(x_1, \ldots ,x_j, \ldots x_n)$$
	is non-constant almost surely; assume for simplicity that all the relevant conditional distributions are equivalent to Lebesgue measure.
	Hence, by Remark \ref{fibers}, these additive noise models will satisfy the converse pairwise Markov property, provided that the true causal graph is a DAG.
\end{remark}

Recall that if $ i \arc j$, then $\epsilon_i$ and $\epsilon_j$ are dependent.  In order to deal with the possibility of arcs, we have the following restriction that allows us to preserve dependence.   We say that $\phi: \mathcal{X}_1 \times \cdots \times \mathcal{X}_n \times \mathcal{E} \to \mathcal{X}_{n+1}$ in an SCM is \emph{noise injective} if for each fixed $x_1, \ldots, x_n$, the function $e \mapsto \phi(x_1, \ldots, x_n,e)$ is injective.

\begin{lemma}
	\label{arc}
	    Consider an SCM, where $j \arc i$.   If both $\phi_i$ and $\phi_j$ are noise injective, then    $i\notci j\cd \an(i,j)$.
	\end{lemma}
\begin{proof}
	Fix $x_{\an(i,j)}$.   Since $j \arc i$, all the parents of $i$ and $i$ are included in $\an(i,j)$ and $\epsilon_i$ and $\epsilon_j$ are dependent.    The conditional distribution of $(X_i, X_j)$ given $X_{\an(i,j)} = x_{\an(i,j)}$ is the law of
		$$ \big( ( \phi_i(x_{\pa(i)}, \epsilon_i),  \phi_j(x_{\pa(j)}, \epsilon_j) \big),$$
		since  the graph is ancestral and the joint law of $(\epsilon_i, \epsilon_j)$ remains the same under this conditioning.
Since $\phi_i$ and $\phi_j$ are noise injective,   it follows that $\phi_i(x_{\pa(i)}, \epsilon_i)$ and    $\phi_j(x_{\pa(j)}, \epsilon_j)$ are dependent if and only if $\epsilon_i$ and $\epsilon_j$ are dependent.
	\end{proof}

From the proof of Lemma \ref{arc}, we see that noise injectivity is a condition which allows us  to preserve the dependencies of the noises; in general if the random variables $\epsilon_i$ and $\epsilon_j$ are dependent, we cannot conclude (like independence) that functions of  $\epsilon_i$ and $\epsilon_j$ will remain dependent.  It is easy to verify that noise injectivity is satisfied by the additive noise models.

\begin{corollary}
	The additive noise models of Remark \ref{additive} satisfy the converse pairwise assumption, without the assumption that  the true causal graph is a DAG.
	\end{corollary}

\begin{proof}	
Immediate from Remark \ref{additive} and Lemma \ref{arc}.
\end{proof}

  We also say that an SCM is \emph{noise injective} if every function in the SCM is noise injective.  Notice that these conditions do not depend on the structure of the causal graph.

 \begin{corollary}
 	\label{pairwise-pairwise}
 	Any  SCM satisfying positivity, with non-constant fibers, and noise injectivity has the converse pairwise Markvov property; furthermore, if the true causal graph is a DAG, then the noise injectivity assumption may be dropped.
 	\end{corollary}
 	
 \begin{proof}
 	Immediate from Proposition \ref{pairwise-distinct}, Lemma \ref{arc}, and  Remark \ref{fibers}.
 	\end{proof}

Hence, under some additional assumptions, the SCMs are not only Markovian, but minimal Markovian to their true causal graph:
\begin{corollary}\label{coro:scm-minmark}
Consider an SCM with positivity, non-constant fibers, noise injectivity,  and  ordered upward- and downward-stability w.r.t.\ $G_0$. Then  $P$ is minimally Markovian to $G_0$; furthermore, if the true causal graph is a DAG, then the noise injectivity assumption may be dropped.
\end{corollary}
\begin{proof}
The proof follows from the Corollary \ref{coro:minmark} together with Theorem \ref{prop:scm-mark} and Corollary \ref{pairwise-pairwise}.
\end{proof}
We can now present our main result for structure learning for SCMs:
\begin{theorem}
	\label{main-scm-cor}
Consider an SCM  with positivity,  non-constant fibers, and a distribution $P$, which satisfies ordered upward- and downward-stability w.r.t.\ the causal graph $G_0$.    Consider natural structure-learning algorithms.
\begin{enumerate}
  \item If $P$ is path-stable, and the SCM is noise injective,  then  $G(P)$, an output of a
  natural
   structure-learning algorithm, is Markov equivalent to $G_0$.
  \item If $G_0$ is a DAG, and $P$ is V-stable, then a DAG-output of
  a natural
  structure-learning algorithm is Markov equivalent to $G_0$.
\end{enumerate}
\end{theorem}
\begin{proof}
The proof follows from Theorem \ref{coro:main}, Theorem \ref{prop:scm-mark}, and Corollary \ref{pairwise-pairwise}.
\end{proof}

In Theorem \ref{main-scm-cor}, the most difficult condition to verify are ordered stabilities, which may not be directly testable.

\subsection{Specialization to specific SCMs}

We consider some more restrictive assumptions on an SCM that will imply the converse pairwise Markov property.

\subsubsection{Discrete random variables}

For discrete distributions, we have the following simple sufficient condition for an SCM to satisfy the converse pairwise Markov property.    We say that an SCM is \emph{discrete} if both the endogenous and exogenous random variables are discrete.

\begin{proposition} Consider an  SCM  that is  discrete  and satisfies positivity.  Suppose that for all the nodes $i$, the cardinality of the support of $\epsilon_i$ is strictly less  than the support of $X_i$.
	\begin{enumerate}
		\item If $P$ is path-stable, and the SCM is noise injective,  then  $G(P)$, an output of
		a natural
		 structure-learning algorithm, is Markov equivalent to $G_0$.
		\item If $G_0$ is a DAG, and $P$ is V-stable, then a DAG-output of a natural  structure-learning algorithm is Markov equivalent to $G_0$.
	\end{enumerate}
\end{proposition}

\begin{proof}From Theorems \ref{coro:main} and  \ref{prop:scm-mark} it suffices to verify the converse pairwise Markov property for a pair of adjacent nodes $i$ and $j$.     Lemma \ref{arc} takes the care of the case where $i \arc j$.  Suppose that $i \fra j$.
	If the support of $\epsilon_i$ is  strictly less  than the support of $X_i$, then positivity implies that \eqref{dist-eq} is not satisfied so that Proposition \ref{pairwise-distinct} gives that the conditional independence statement $i\ci j\cd \an(i,j)$ does not hold.
\end{proof}

The assumption on the cardinality of the supports has the interpretation that  information at a node is {\em not} mostly noise.

In the next lemma, we show that under some assumptions, the failure of the converse pairwise Markov property for $j \fra i$, forces the uniform distribution on both $X_i$ and $\epsilon_i$.
Let $\phi: \mathcal{X}_1 \times \cdots \times \mathcal{X}_n \times \mathcal{E} \to \mathcal{X}_{n+1}$   be a noise injective function in a discrete SCM.
Given $x_1, \ldots, x_j, \ldots, x_n$, let  $\Phi_{x_j} ^{x_1, \ldots, x_{j-1}, x_{j+1}, \ldots, x_n} = \Phi_{x_j} ^{x^j}$ denote the inverse of $\phi$.  Suppose that for all choices of $1 \leq j \leq n$ and  $x^j$ there exists $x_{n+1}^* \in \mathcal{X}_{n+1}$ such that
$$\{  \Phi_{x_j} ^{x^j} (x_{n+1} ^*) : x_{j} \in \mathcal{X}_{j} \} = \mathcal{E}.$$  Then we say that $\phi$ is \emph{noise surjective}.

\begin{proposition}
	\label{uniform-discrete}
	Consider the following SCM, where $j \fra i$ and $$X_i= \phi_i(X_{\pa(i) \setminus \{j\} }, X_j, \epsilon_i).$$  Suppose that the joint distribution   is discrete and satisfies positivity.   Suppose also that $\phi_i$ is noise injective and surjective.
	If  $i\ci j\cd \an(i,j)$, then  $\epsilon_i$ and $X_j$ are both uniformly distributed and $\mathcal{E}_i$ and $\mathcal{X}_i$ have the same cardinality.
\end{proposition}

\begin{proof}  Fix $x_{\an(i,j) }$.  Note that all the parents of $i$, save $j$,  are included in $\an(i,j)$,
	By Proposition	\ref{pairwise-distinct}, we have
	\begin{eqnarray*}
	\mathbb{P}(X_i =x_i| X_{\an(i,j) }  =  x_{\an(i,j) } )
	 &=& \mathbb{P} (\phi_i(x_{{\pa(i) \setminus \{j\} })},x_j,\epsilon_i) =x_i) \\
	 &=& \mathbb{P} (\epsilon_i =\Phi_{x_j} ^{x_{\pa(i) \setminus \{j\} }}  (x_i) )
	 \end{eqnarray*}
	 for all $x_i \in \mathcal{X}_i$ and all $x_j \in \mathcal{X}_j$,  by positivity.
By noise surjectivity, substituting  $x_i = x_i^*$, we deduce that  $\epsilon_i$ is uniformly distributed, with finite support $\mathcal{E}_i$. Thus
$$	\mathbb{P}(X_i =x_i| X_{\an(i,j) }  =  x_{\an(i,j) } ) = \frac{1}{ \# \mathcal{E}_i}$$
for all $x_i$.  Hence it
 follows that $X_i$ is also uniformly distributed.
	\end{proof}

\begin{corollary}
\label{unif-noisea}
 Consider a discrete SCM satisfying positivity.  Suppose that all the functions are noise injective and noise surjective.  If each of the  exogenous random variables are not uniformly distributed, then the converse pairwise Markov property is satisfied.
	\end{corollary}

\begin{proof}
	Immediate from Proposition \ref{uniform-discrete} and Lemma \ref{arc}.
	\end{proof}

\begin{corollary}
\label{unif-noiseb}
 Consider an  SCM  that is  discrete  and satisfies positivity. Suppose that all the functions are noise injective and noise surjective.  Assume  each of the  exogenous random variables are not uniformly distributed.
	\begin{enumerate}
		\item If $P$ is path-stable,   then  $G(P)$, an output of
		a natural
		 structure-learning algorithm, is Markov equivalent to $G_0$.
		\item If $G_0$ is a DAG, and $P$ is V-stable, then a DAG-output of a natural  structure-learning algorithm is Markov equivalent to $G_0$.
	\end{enumerate}
\end{corollary}

\begin{proof}
The previous corollary assures us that the converse pairwise Markov property is satisfied, so that the result follows from Theorems \ref{coro:main} and \ref{prop:scm-mark}.
\end{proof}

Returning to Example \ref{mod}, it is easy to verify that the simple SCM is both noise injective and surjective.   Notice that if $\epsilon_1$ is chosen with Bernoulli parameter $p \not =\tfrac{1}{2}$, and all other choices are left the same, then $X_2$ is no longer independent of $X_1$.   Corollaries \ref{unif-noisea} and \ref{unif-noiseb} suggest that  if it is possible to avoid or exclude the possibility that the exogenous random variables are  uniformly distributed, then a natural structure-learning algorithm is more likely to work.

\subsubsection{Continuous random variables}

We note that the modular addition in Example \ref{mod} is not continuous.  We show that under some regularity conditions, up to a change of coordinates, only trivial functions can result in independence.
Given $\phi: \mathcal{X}_1 \times \cdots \times \mathcal{X}_n \times \mathcal{E} \to \mathcal{X}_{n+1}$, 	we say that $\psi$ is a \emph{$j$-change of coordinates} if there exists $g: [-1/2, 1/2]  \to \mathcal{E}$ such that  for all $x^j \in \mathcal{X}^j$, there exists $h_{x^j} : \mathcal{X}_{n+1} \to [-1/2, 1/2]$ such that if  $U$ is uniformly distributed on $[-1/2, 1/2]$ and
$$\psi(x_j, u) :=  h_{x^j} \circ \phi(x_1, \ldots, x_j, \ldots, x_n, g(u)),$$
then $g(U)$ has the law of $\epsilon$ and  $\psi(x_j, \epsilon)$ is uniformly distributed on $[-1/2, 1/2]$ for all $x_j \in \mathcal{X}_j$.  We say that $\psi$ is \emph{regular} if for fixed $x_j$, the function $\zeta(u) =  \psi(x_j, u)$ is injective and $\zeta(u)$ is  monotone; note that under these conditions an easy argument given in Lemma \ref{uniform-regular} shows that  $\psi$ is a \emph{reflection} of  the form $\psi(x_j, u) = u$ or $\psi(x_j, u) =- u$ for all $x_j \in \mathcal{X}_j$ and all $u \in [-1/2, 1/2]$.   We say that $\phi$ is \emph{non-trivial} if for each $j$, there does not exist a regular $j$-change of coordinates.
We say that an SCM is  \emph{non-trivial} if every function in the SCM is non-trivial.
We say that a function $\phi$ in an SCM is \emph{noise monotone} if for all $x_1, \ldots, x_n$, the function  $e \mapsto \phi(x_1, \ldots, x_n,e) \in \mathcal{X}_{n+1} \subset \mathbb{R}$ is monotone in $e \in \mathcal{E} \subset \mathbb{R}$.

\begin{proposition}
	\label{uniform-regular}
	Consider the following SCM, where $j \fra i$ and $$X_i= \phi_i(X_{\pa(i) \setminus \{j\} }, X_j, \epsilon_i).$$  Suppose that $\phi_i$ is noise injective and noise monotone.    Also assume that $\epsilon_i$ is a continuous random variable, and that for almost all $x_{\an(i,j)}$,  the conditional distribution law  of $X_i$ given $X_{\an(i,j)}$ is continuous with respect to Lebesgue measure.  If  $i\ci j\cd \an(i,j)$, then   $\phi_i$ admits a $j$-regular change of coordinates and any such change of coordinates  is a  reflection.
\end{proposition}

\begin{proof}
	Recall that if $F$ is the cumulative distribution function  for a continuous random variable $Z$, then  $F$ is non-decreasing and $F(Z) - \tfrac{1}{2}$ is uniformly distributed on $[-1/2, 1/2]$.  In addition, the inverse transform method tells us that $F^{-1}(U+\tfrac{1}{2})$ has law $Z$ if $U$ is uniformly distributed on $[-1/2,1/2]$.
	Hence by Proposition \ref{pairwise-distinct}, and the absolute continuity assumptions, it follows that $\phi_i$ has a $j$-regular change of coordinates $\psi$.    Fix $x_j \in \mathcal{X}_j$ and write $\psi(x_j, u) = \psi(u)$.      If $\psi$ is increasing in $u$, then for $z \in [-1/2, 1/2]$, we have for $U$ uniformly distributed on $[-1/2, 1/2]$ that
	$$ z+1/2= \mathbb{P}(\psi(U) \leq z) = \mathbb{P}( U \leq \psi^{-1}(z)) = \psi^{-1}(z) + 1/2,$$
	from which it follows that $\psi(x_j,z) = z$.
	If $\psi$ is decreasing,  then a similar argument gives that $\psi(x_j,z) = -z$.
\end{proof}

\begin{corollary}
If the law of an SCM is given by a jointly continuous distribution, and the SCM is noise injective and monotone, and is non-trivial, then it satisfies the converse pairwise Markov property.  Moreover,
\begin{enumerate}
		\item If $P$ is path-stable,   then  $G(P)$, an output of
		a natural
		 structure-learning algorithm, is Markov equivalent to $G_0$.
		\item If $G_0$ is a DAG, and $P$ is V-stable, then a DAG-output of a natural  structure-learning algorithm is Markov equivalent to $G_0$.
	\end{enumerate}

	\end{corollary}

\begin{proof}
	The converse pairwise Markov property is immediate from Proposition \ref{uniform-regular} and Lemma \ref{arc}.  Hence the latter assertions that a natural structure-learning algorithm works follow from Theorems \ref{coro:main} and \ref{prop:scm-mark}.
	\end{proof}

\section{Summary and discussion}\label{sec:summary}
We formalized the problem of structure learning from independence statements to find a graph that is Markov equivalent to the ``true causal graph," which is the unknown graph to which the distribution is assumed to be Markovian. In fact, all the results presented here are also correct for an abstract independence model $\mathcal{J}$ instead of a probability distribution $P$.

A main message of this paper is that an important property of constraint-based structure-learning algorithms is to direct the edges of the skeleton of the probability distribution $P$ so that ordered upward- and downward-stability are satisfied by $P$ w.r.t.\ the generated graph. We call such algorithms a natural constraint-based structure-learning algorithm. Among the conditions that constitute the widely-used assumption of faithfulness, provided in Proposition \ref{thm:2}, ordered upward- and downward-stability are the only ones that induce an ordering, which is what is used for directing the edges in structure learning. Since, under faithfulness,  ordered upward- and downward-stability are necessary for the output of the algorithm, all exact structure-learning algorithms that assume faithfulness in the literature are members of the class of natural structure-learning algorithms. Therefore, for any (newly) proposed algorithm, to prove that the algorithm works, it is sufficient to test whether $P$ satisfies  ordered upward- and downward-stability w.r.t.\ its output.

Following this line of thinking, this result could lead to finding the ``optimal" natural constraint-based structure-learning algorithm. The idea is that such an algorithm would be designed to solely direct the edges of the skeleton to ensure an ordering w.r.t.\ which  ordered upward- and downward-stability are satisfied.  Hence, for certain recorded independence statements, one should perform some additional tests on whether independence is preserved by  certain additions to the conditioning set of those statements. The results of such tests determine the direction of the edges of the output. This is a direction for future work.


Moreover, we provided two assumptions, both of which are satisfied by faithfulness, in order to ensure that the skeleton of the output of a natural structure-learning algorithm is the same as the skeleton of the causal graph. One condition is the converse pairwise Markov property and the other one is satisfied by the assumption that the distribution satisfies ordered upward- and downward-stability w.r.t.\ the causal graph. In order to ensure that directing the edges of the skeleton of the output will result to a Markov equivalent graph to the causal graph, a sufficient condition is path-stability in the general case of directed ancestral graphs, and V-stability in the case where we only deal with DAGs.

In contrast to well-known conditions in the literature such as orientation-faithfulness, path- and V-stability are only  properties of the distribution and \emph{not} the causal graph. V-stability is also weaker than singleton-transitivity, and consequently faithfulness. However, there are cases where faithfulness is satisfied, but not path-stability. Path- and V-stability are easier to test than singleton-transitivity as one needs to only test independence for the endpoints of specific paths or V-configurations in $\sk(P)$ respectively. In addition to singleton-transitivity, out of the equivalent conditions to faithfulness, provided in \citet{sad17}, intersection and composition properties are not needed. In the case of a Gaussian distribution, however, all these assumptions are automatically satisfied, and therefore, our assumptions imply faithfulness.

It is commonplace to perform structure learning for SCMs. We have specialized the theory for SCMs. For SCMs that we associate with directed ancestral graphs, we showed that the main assumption of $P$ being Markovian to the causal graph is always satisfied.


In addition, we provided simple sufficient conditions on the SCM that imply the converse pairwise Markov property. Although converse pairwise Markov property is indeed related to the causal graph, these  sufficient conditions are not. Our key condition of having non-constant fibers has a simple interpretation when the  function of an SCM has only one discrete endogenous variable and one exogenous variable: there is a choice of two values for the  endogenous variable, whereby the distributions of the random variables given by evaluating the function at these two fixed values and the exogenous random variable are distinct.      In the general case where the SCM may not be a DAG, we need to impose a further strong injectivity condition on the functions.

The assumption of distributions (of the SCM) satisfying ordered upward- and downward-stability w.r.t.\ the causal graph is, in principle, not directly testable. However, in our primary study, we observe that it is  difficult to come up with natural  examples of distributions that do not satisfy these conditions when dealing with SCMs. Nevertheless, a thorough study of these conditions for SCMs and other causal models is needed, and is subject of future work.

\acks{We are grateful to Steffen Lauritzen for providing Example \ref{ex:1}, to Joris Mooij for stimulating discussions with the first author especially on the examples involving modular arithmetic, and to Jonas Peters for hosting a first author's visit, during which some of the ideas used in this paper were discussed. We  are also truly thankful to the associate editor and three anonymous reviewers for detailed and helpful comments.}






\vskip 0.2in
\bibliography{bib}

\end{document}